\documentclass[11pt]{article}
\usepackage{amssymb}
\usepackage{times}

\title{How can we recognize potentially $\bormxi$ subsets of the plane?\indent}
\author{Dominique LECOMTE}
\date{June 2006}

\def\ufootnote#1{\let\savedthfn\thefootnote\let\thefootnote\relax
\footnote{#1}\let\thefootnote\savedthfn\addtocounter{footnote}{-1}}

\newcommand{\Ana}{{\it\Sigma}^{1}_{1}}
\newcommand{\Ca}{{\it\Pi}^{1}_{1}}
\newcommand{\Borel}{{\it\Delta}^{1}_{1}}

\newcommand{\borel}{{\bf\Delta}^{1}_{1}}
\newcommand{\boraone}{{\bf\Sigma}^{0}_{1}}

\newcommand{\boraxi}{{\bf\Sigma}^{0}_{\xi}}

\newcommand{\bormz}{{\bf\Pi}^{0}_{0}}
\newcommand{\bormone}{{\bf\Pi}^{0}_{1}}
\newcommand{\bormtwo}{{\bf\Pi}^{0}_{2}}

\newcommand{\borapx}{{\bf\Sigma}^{0}_{1+\xi}}

\newcommand{\bormpx}{{\bf\Pi}^{0}_{1+\xi}}
\newcommand{\bormlxi}{{\bf\Pi}^{0}_{<\xi}}
\newcommand{\bormlpxi}{{\bf\Pi}^{0}_{<1+\xi}}

\newcommand{\bormep}{{\bf\Pi}^{0}_{\eta +1}}
\newcommand{\bormxi}{{\bf\Pi}^{0}_{\xi}}
\newcommand{\borme}{{\bf\Pi}^{0}_{\eta}}

\newtheorem{thm} {Theorem} [section]
\newtheorem{defi} [thm] {Definition}
\newtheorem{cor} [thm] {Corollary}
\newtheorem{lem} [thm] {Lemma}
\newtheorem{prop} [thm] {Proposition}
\newtheorem{them} {Theorem} [subsection]
\newtheorem{defin} [them] {Definition}

\newtheorem{lemm} [them] {Lemma}

\begin{document}

\maketitle

\centerline{$\bullet$ Universit\' e Paris 6, Institut de Math\'ematiques de Jussieu, tour 46-0, bo\^\i te 186,}

\centerline{4, place Jussieu, 75 252 Paris Cedex 05, France.}

\centerline{dominique.lecomte@upmc.fr}\bigskip

\centerline{$\bullet$ Universit\'e de Picardie, I.U.T. de l'Oise, site de Creil,}

\centerline{13, all\'ee de la fa\"\i encerie, 60 107 Creil, France.}\bigskip\bigskip\bigskip\bigskip\bigskip\bigskip\bigskip

\ufootnote{{\it 2000 Mathematics Subject Classification.}~Primary: 03E15, Secondary: 54H05}

\ufootnote{{\it Keywords and phrases.}~Potentially, Baire classes, Reduction, Hurewicz's Theorem}

\ufootnote{{\bf Acknowledgements.}~I would like to thank G. Debs for his nice 
contributions to this paper. The main results proved in this paper are 
the solution to the main question (partially solved) in my thesis. 
This question was asked by A. Louveau fifteen years ago. I would like to 
thank him very much for that. I also thank them for their remarks during the 
lectures I gave in the descriptive set theory seminar of the University Paris 6, 
where I proved some of the results in this paper. These remarks improved the quality of this article.}

\noindent {\bf Abstract.} Let $\xi\!\geq\! 1$ be a countable ordinal. We study the Borel 
subsets of the plane that can be made $\bormxi$ by refining the Polish topology 
on the real line. These sets are called potentially $\bormxi$. We give a 
Hurewicz-like test to recognize potentially $\bormxi$ sets.

\vfill\eject

\section{$\!\!\!\!\!\!$ Introduction.}\indent

 The reader should see [K] for the descriptive set theoretic notation used in 
this paper. This work is the continuation of a study started in [L1]-[L5], and is announced in [L6]. 
The usual notion of comparison for Borel equivalence relations $E\!\subseteq\! X^2$ and 
${E'\!\subseteq\! {X'}^2}$ on Polish spaces is the Borel reducibility 
quasi-order:
$$E\leq_B E'~\Leftrightarrow ~\exists u\!:\! X\!\rightarrow\! X'~
\hbox{\rm Borel~with}~E\! =\! (u\!\times\! u)^{-1}(E')$$
(recall that a quasi-order is a reflexive and transitive relation). Note that 
this makes sense even if $E$, $E'$ are not equivalence relations.
It is known that if $(B_n)$ is a sequence of Borel subsets of $X$, then there is 
a finer Polish topology on $X$ making the $B_n$'s clopen (see exercise 13.5 in 
[K]). So assume that $E\leq_B E'$, and let $\sigma$ be a finer Polish topology on 
$X$ making $u$ continuous. If moreover $E'$ is in some Baire class $\Gamma$, 
then $E\!\in\!\Gamma ([X,\sigma ]^2)$. This motivates the following (see [Lo2]):

\begin{defi} (Louveau) Let $X$, $Y$ be Polish spaces, $A$ a Borel 
subset of $X\!\times\! Y$, and $\Gamma$ a Baire (or Wadge) class. We say that $A$ is 
$\underline{potentially~in~\Gamma}$ 
$\big($denoted ${A\!\in\!\hbox{\it pot}(\Gamma)\big)}$ iff there is a finer 
Polish topology $\sigma$ $($resp., $\tau )$ on $X$ $($resp., $Y)$ with  
$A\!\in\!\Gamma ([X,\sigma ]\!\times\! [Y,\tau ])$.\end{defi}

 This notion is a natural invariant for $\leq_B$: if 
$E'$ is $\hbox{\rm pot}(\Gamma)$ and $E\leq_B E'$, then 
$E$ is $\hbox{\rm pot}(\Gamma)$. Using this notion, A. Louveau proved that the 
collection of $\boraxi$ equivalence 
relations is not cofinal for $\leq_B$, and deduces from this the non 
existence of a maximum Borel equivalence relation for $\leq_B$ (this non existence 
result is due to H. Friedman and L. Stanley). More 
recently, G. Hjorth, A. Kechris and A. Louveau determined the 
potential classes of the Borel equivalence relations induced by Borel 
actions of closed subgroups of the symmetric group (see [Hj-K-Lo]).\bigskip

 A standard way to see that a set is complicated is to note that it is more 
complicated than a well-known example. For instance, we have the following (see 
[SR]):

\begin{thm} (Hurewicz) Let $P_{f}\! :=\!\{\alpha\!\in\! 2^\omega\mid
\exists n\!\in\!\omega ~\forall m\geq n\ \ \alpha (m)\! =\! 0\}$, $X$ be 
a Polish space, and $A$ a Borel subset of $X$. Then exactly one of the 
following holds:

\noindent (a) The set $A$ is $\bormtwo(X)$.

\noindent (b) There is $u\! :\! 2^\omega\!\rightarrow\! X$ continuous and one-to-one with 
$P_{f}\! =\! u^{-1}(A)$.\end{thm}

 This result has been generalized to all Baire classes (see [Lo-SR]). We state 
this generalization in two parts:

\begin{thm} (Louveau-Saint Raymond) Let $\xi\! <\!\omega_{1}$, 
$S\!\in\!\borapx (2^{\omega} )$, $X$ be a Polish space, and $A$, $B$ 
disjoint analytic subsets of $X$. Then one of the following holds:

\noindent (a) The set $A$ is separable from $B$ by a $\bormpx (X)$ set.

\noindent (b) There is $u\! :\! 2^{\omega}\!\rightarrow\! X$ continuous 
with $S\!\subseteq\! u^{-1}(A)$ and $2^{\omega}\!\setminus\! S\!\subseteq\! u^{-1}(B)$.

 If we moreover assume that $S\!\notin\!\bormpx$, then this is a dichotomy.\end{thm}

 Note that in this dichotomy, we can have $u$ one-to-one if $\xi\!\geq\! 2$. This is not possible if 
$\xi\! <\! 2$.

\begin{thm} There is a concrete example of a set 
$S_{1+\xi}\!\in\!\borapx (2^{\omega} )\!\setminus\!\bormpx (2^{\omega} )$, for each 
$\xi\! <\!\omega_{1}$.\end{thm}
 
 We try to adapt these results to the Borel subsets of the plane. 
 
\vfill\eject
 
 The following result is proved in [H-K-Lo]:

\begin{thm} (Harrington-Kechris-Louveau) Let $X$ be a Polish space, $E$ a 
Borel equivalence relation on $X$, and $E_{0}\! :=\!\{(\alpha ,\beta )
\!\in\! 2^{\omega}\!\times\! 2^{\omega}\mid\exists n\!\in\!{\omega}\ \ 
\forall m\!\geq\! n\ \ \alpha (m)\! =\!\beta (m)\}$. Then exactly one of the 
following holds:

\noindent (a) The relation $E$ is $\hbox{\it pot}(\bormone)$.

\noindent (b) $E_{0}\leq_B E$ (with $u$ continuous and one-to-one).\end{thm}

 For the Borel subsets of the plane, we need some other notions of 
comparison. Let $X$, $Y$, $X'$, $Y'$ be Polish spaces, and $A$ (resp., 
$A'$) a Borel subset of $X\!\times\! Y$ (resp., $X'\!\times\! Y'$). We set
$$A\leq^r_{B}A'~\Leftrightarrow ~\exists u\! :\! X\!\rightarrow\! X'\ 
\ \exists v\! :\! Y\!\rightarrow\! Y'\ \ \hbox{\rm Borel\ with}\ A\! =\! 
(u\!\times\! v)^{-1}(A').$$
The following result is proved in [L1]:
 
\begin{thm} Let $\Delta (2^{\omega} )\! :=\!\{(\alpha ,\beta )\!\in
\! 2^{\omega}\!\times\! 2^{\omega}\mid\alpha\! =\!\beta\}$,  
${L_0\! :=\!\{(\alpha ,\beta )\!\in\! 2^{\omega}\!\times\! 2^{\omega}\mid
\alpha\! <_{\hbox{\it lex}}\!\beta\}}$, $X$, $Y$ be Polish spaces, and $A$ a 
$\hbox{\it pot}\big(\check D_2(\boraone)\big)$ subset of $X\!\times\! Y$. 
Then exactly one of the following holds:

\noindent (a) The set $A$ is $\hbox{\it pot}(\bormone)$.

\noindent (b) $\neg\Delta (2^{\omega} )\leq^r_{B} A$ or $L_{0}\leq^r_{B} A$ (with 
$u$, $v$ continuous and one-to-one).\end{thm}
 
 The class $\check D_2(\boraone)$ is the class of unions of a closed set and of 
an open set. Things become more complicated at the level 
$D_{2}(\boraone)$ of differences of two open sets (see [L5]):

\begin{thm} (a) There is a perfect $\leq^r_{B}$-antichain 
$(A_{\alpha})_{\alpha\in 2^{\omega}}\!\subseteq\! 
D_{2}(\boraone )(2^{\omega}\!\times\! 2^{\omega})$ 
such that $A_{\alpha}$ is $\leq^r_{B}$-minimal among 
$\borel\!\setminus\!\hbox{\it pot}(\bormone )$ sets, for any 
$\alpha\!\in\! 2^{\omega}$.

\noindent (b) There is a perfect $\leq_{B}$-antichain 
$(R_{\alpha})_{\alpha\in 2^{\omega}}$ such that $R_{\alpha}$ is 
$\leq_{B}$-minimal among 
$\borel\!\setminus\! \hbox{\it pot}(\bormone)$ sets, for any 
$\alpha\!\in\! 2^{\omega}$. Moreover, $(R_{\alpha})_{\alpha\in 2^{\omega}}$ 
can be taken to be a subclass of any of the following classes:

- Graphs  (i.e., irreflexive and symmetric relations).

- Oriented graphs  (i.e., irreflexive and antisymmetric relations).

- Quasi-orders.

- Partial orders (i.e., reflexive, antisymmetric and transitive 
relations).\end{thm}
 
  In other words, the case of equivalence relations, for which we have a 
unique (up to bi-reducibili-ty) minimal non potentially closed element with Theorem 1.5, 
is very specific. Theorem 1.7.(b) says, among other things, 
that the mixture between symmetry and transitivity is very strong. Theorem 1.7.(a) 
shows that the classical notions of reduction (on the whole product) don't work, 
at least at the first level. So we must find another notion of comparison. The 
following result is proved in [L5]:

\begin{thm} There is $S_{1}\!\in\!\borel (2^\omega\!\times\! 2^\omega )$ such that 
for any Polish spaces $X$, $Y$, and for any Borel subset $A$ of $X\!\times\! Y$, 
exactly one of the following holds:

\noindent (a) The set $A$ is $\hbox{\it pot}(\bormone )$.

\noindent (b) There are $u\! :\! 2^\omega\!\rightarrow\! X$ and 
$v\! :\! 2^\omega\!\rightarrow\! Y$ continuous satisfying the inclusions 
$S_{1}\!\subseteq\! (u\!\times\! v)^{-1}(A)$ and 
$\overline{S_{1}}\!\setminus\! S_1\!\subseteq\! (u\!\times\! v)^{-1}(\neg A)$.

 Moreover, we can neither replace $\overline{S_{1}}\!\setminus\! S_{1}$ with $\neg S_{1}$, 
nor ensure that $u$ and $v$ are one-to-one.\end{thm}

\vfill\eject

 So we get a minimum non-potentially closed set if we do not ask for 
a reduction on the whole product. We will show that this dichotomy is true for 
each countable ordinal $\xi\!\geq\! 1$. The result is actually stronger than 
that. First the $A_{\xi}$'s are concrete examples. Secondly it is better to state 
that the reduction in condition (b) holds in the set $\lceil T\rceil$ of the 
branches of some tree $T$ that does not depend on $\xi$, rather than 
$\overline{A_{\xi}}$. Finally, to get the full strength of the result, it is 
better to split it in two parts. We need some notation and a definition:\bigskip

\noindent\bf Notation.\rm\ If ${\cal F}_{0}$, ${\cal F}_{1}$ are finite sets and 
${\cal T}\!\subseteq\! {\cal F}_{0}\!\times\! {\cal F}_{1}$, we denote by 
$G_{\cal T}$ the bipartite graph with set of vertices the sum 
${\cal F}_{0}\!\oplus\! {\cal F}_{1}$, and with set of edges 
$$\Big\{\{ (f_{0},0) , (f_{1},1)\}\!\subseteq {\cal F}_{0}\!\oplus\! {\cal F}_{1}
\mid (f_{0},f_{1})\!\in\! {\cal T}\Big\}.$$
(see [B] for basic notions about graphs). In the sequel, we will denote 
$\overline{f_\varepsilon}\! :=\! (f_\varepsilon ,\varepsilon )$.

\begin{defi} We say that a tree $T$ on $2\!\times\! 2$ is a 
$\underline{tree~with~acyclic~levels}$ if, for each integer $p$, the graph 
$G_{{\cal T}_{p}}$, associated with 
${\cal T}_{p}\! :=\! T\cap (2^p\!\times\! 2^p)\!\subseteq\! 2^p\!\times\! 2^p$, 
is acyclic.\end{defi}

 Now we can state the main results proved in this paper:

\begin{thm} (Debs-Lecomte) Let $T$ be a tree with acyclic  levels, 
$\xi\! <\!\omega_{1}$, $S\!\in\!\borapx (\lceil T\rceil)$, $X$, $Y$ Polish 
spaces, and $A$, $B$ disjoint analytic subsets of $X\!\times\! Y$. Then one of 
the following holds:

\noindent (a) The set $A$ is separable from $B$ by a $\hbox{\it pot}(\bormpx )$ set.

\noindent (b) There are $u\! :\! 2^\omega\!\rightarrow\! X$ and 
$v\! :\! 2^\omega\!\rightarrow\! Y$ 
continuous with $S\!\subseteq\! (u\!\times\! v)^{-1}(A)$ and 
$\lceil T\rceil\!\setminus\! S\!\subseteq\! (u\!\times\! v)^{-1}(B)$.

 If we moreover assume that $S\!\notin\! \hbox{\it pot}(\bormpx )$, then this is a dichotomy.\end{thm}

 Note that we can deduce Theorem 1.3 from the proof of Theorem 1.10. Theorem 1.10 is the ana-logous of Theorem 1.3 in dimension two. The proofs of Theorem 1.3 in [Lo-SR], and also Theorem III-2.1 in [D-SR], use games. This is not the case here, so that we get a new proof of Theorem 1.3. 

\begin{thm} We can find concrete examples of a tree $T$ with acyclic levels, together with sets 
$S_{1+\xi}\!\in\!\borapx (\lceil T\rceil)\!\setminus\!\hbox{\it pot}(\bormpx )$, for each 
$\xi\! <\!\omega_{1}$.\end{thm}

 The following corollary has initially been shown by D. Lecomte when $1\! +\!\xi$ is a 
successor ordinal. Then G. Debs showed it when $1\! +\!\xi$ is a limit ordinal. 

\begin{cor} (Debs-Lecomte) Let $\xi\! <\!\omega_1$. There is 
$S\!\in\!\borel (2^\omega\!\times\! 2^\omega )$ such that for any 
Polish spaces $X$, $Y$, and for any disjoint analytic subsets $A$, $B$ of 
$X\!\times\! Y$, exactly one of the following holds:

\noindent (a) The set $A$ is separable from $B$ by a $\hbox{\it pot}(\bormpx )$ set.

\noindent (b) There are $u\! :\! 2^\omega\!\rightarrow\! X$ and 
$v\! :\! 2^\omega\!\rightarrow\! Y$ continuous with 
$S\!\subseteq\! (u\!\times\! v)^{-1}(A)$ and 
$\overline{S}\!\setminus\! S\!\subseteq\! (u\!\times\! v)^{-1}(B)$.\end{cor}
 
 Theorem 1.8 shows that we cannot replace $\overline{S}\!\setminus\! S$ with $\neg S$ in Corollary 1.12 when $\xi\! =\! 0$. G. Debs found a simpler proof, which moreover works in the general case:
 
\begin{thm} (Debs) We cannot replace $\overline{S}\!\setminus\! S$ with $\neg S$ in Corollary 1.12.
\end{thm}
 
 Once again, some cycles are involved, so that the acyclicity is essentially necessary and sufficient in Corollary 1.12 (even if we have two different notions of acyclicity). G. Debs proved very recently that we 
can have $u$ and $v$ one-to-one in Corollary 1.12 if $\xi\!\geq\! 2$. This is not possible if $\xi\! <\! 2$ (see Theorem 1.8 when $\xi\! =\! 0$, and Theorem 15 in [L4] when $\xi\! =\! 1$). 
 
\vfill\eject

 This paper is organized as follows:\bigskip
 
\noindent - In Section 2 we recall the material used to state the 
representation theorem of Borel sets proved in [D-SR]. We use it to 
prove Theorem 1.10, also in this section. To do this we assume some 
results proved in [Lo2]. We also prove Theorem 1.13.\bigskip

\noindent - In Section 3 we prove Theorem 1.11.\bigskip

\noindent - We use some tools of effective descriptive set theory 
(the reader should see [M] for the basic notions about it). In Section 4 we give an 
alternative proof of the results in [Lo2] that we assumed in Section 2. This 
leads to the following:

\begin{thm} (Debs-Lecomte-Louveau) Let $T$ given by Theorem 1.11, 
$\xi\! <\!\omega^{\hbox{\it CK}}_{1}$, $S$ given by Theorem 1.11, $X$, $Y$ 
be recursively presented Polish spaces, and $A$, $B$ disjoint $\Ana$ subsets of 
$X\!\times\! Y$. Then the following are equivalent:

\noindent (a) The set $A$ cannot be separated from $B$ by a $\hbox{\it pot}(\bormpx )$ set.

\noindent (b) The set $A$ cannot be separated from $B$ by a 
$\Borel\cap\hbox{\it pot}(\bormpx )$ set.

\noindent (c) There are $u\! :\! 2^\omega\!\rightarrow\! X$ and 
$v\! :\! 2^\omega\!\rightarrow\! Y$ 
continuous with $S\!\subseteq\! (u\!\times\! v)^{-1}(A)$ and 
$\lceil T\rceil\!\setminus\! S\!\subseteq\! (u\!\times\! v)^{-1}(B)$.
\end{thm}

 The equivalence between (a) and (b) is proved in [Lo2]. We will actually prove 
more than Theo-rem 1.14, with some additional notation that will be introduced later. 
Among other things, we will use the fact that the set of codes for $\Borel$ and 
$\hbox{\rm pot}(\bormpx )$ sets is $\Ca$.

\section{$\!\!\!\!\!\!$ Proof of Theorem 1.10.}

\subsection{$\!\!\!\!\!\!$ Acyclicity.}\indent

 In this subsection we prove a result that will be used later to show Theorem 
1.10. This is the place where the essence of the notion of a tree with acyclic 
levels is really used. We will also prove that we cannot have a reduction on the 
whole product, using some cycles. Some of the arguments used in the initial proof 
of Corollary 1.12 by D. Lecomte (when $1\! +\!\xi$ is a successor ordinal) are 
replaced here by Lemma 2.1.2 below.

\begin{defin} (Debs) Let ${\cal F}_{0}$, ${\cal F}_{1}$, $X_{0}$, $X_{1}$ be sets, 
${\cal T}\!\subseteq\! {\cal F}_{0}\!\times\! {\cal F}_{1}$ and 
$\Psi\! :\! {\cal F}_{0}\!\times\! {\cal F}_{1}\!\rightarrow\! 
2^{X_{0}\times X_{1}}$. We say that 
$\psi\! =\!\psi_{0}\!\times\!\psi_{1}\! :\! {\cal F}_{0}\!\times\! {\cal F}_{1}\!
\rightarrow\! X_{0}\!\times\! X_{1}$ is a 
$\underline{\pi\! -\! selector~on~{\cal T}~for~\Psi}$ if:

\noindent (a) $\psi (f_{0},f_{1})\! =\! [\psi_{0} (f_{0}),\psi_{1} (f_{1})]$, for each 
$(f_{0},f_{1})\!\in\! {\cal F}_{0}\!\times\! {\cal F}_{1}$.

\noindent (b) $\psi (t)\!\in\!\Psi (t)$, for each $t\!\in\! {\cal T}$.\end{defin}

\noindent\bf Notation.\rm ~Let $X$ be a recursively presented Polish 
space. We denote by ${\it\Delta}_{X}$ the topology on $X$ generated by $\Borel 
(X)$. This topology is Polish (see (iii)$\Rightarrow$(i) in the proof of Theorem 
3.4 in [Lo2]). We set $\tau_{1}\!:=\! {\it\Delta}_{X}\!\times\! {\it\Delta}_{Y}$ 
if $Y$ is also a recursively presented Polish space.

\begin{lemm} (Debs) Let ${\cal F}_{0}$, ${\cal F}_{1}$ be finite sets, 
${\cal T}\!\subseteq\! {\cal F}_{0}\!\times\! {\cal F}_{1}$ such that the 
graph $G_{\cal T}$ associated with $\cal T$ is acyclic, $X_{0}$, $X_{1}$ 
recursively presented Polish spaces, 
$\Psi\! :\! {\cal F}_{0}\!\times\! {\cal F}_{1}\!\rightarrow\!\Ana 
(X_{0}\!\times\! X_{1})$, and 
$\overline{\Psi}\! :\! {\cal F}_{0}\!\times\! {\cal F}_{1}\!\rightarrow\!\Ana 
(X_{0}\!\times\! X_{1})$ defined by 
$\overline{\Psi}(t)\! :=\!\overline{\Psi (t)}^{\tau_{1}}$. Then $\Psi$ 
admits a $\pi$-selector on $\cal T$ if $\overline{\Psi}$ does.\end{lemm}

\vfill\eject

\noindent\bf Proof.\rm ~(a) Let $t_{0}\! :=\! (f_{0},f_{1})\!\in\! {\cal T}$, and 
$\Phi\! :\! {\cal F}_{0}\!\times\! {\cal F}_{1}\!\rightarrow\!\Ana (X_{0}\!\times
\! X_{1})$. We assume that $\Phi (t)\! =\!\Psi (t)$ if $t\!\not=\! t_{0}$, and 
that $\Phi (t_{0})\!\subseteq\!\overline{\Psi (t_{0})}^{\tau_{1}}$. We first prove 
that $\Psi$ admits a $\pi$-selector on $\cal T$ if $\Phi$ does.\bigskip

\noindent $\bullet$ Fix a $\pi$-selector $\tilde\varphi$ on $\cal T$ for $\Phi$. 
We define $\Ana$ sets $U_{\varepsilon}$, for $\varepsilon\!\in\! 2$, by
$$U_{\varepsilon} :=\{\ x\!\in\! X_{\varepsilon}\mid\exists\varphi\! :\! 
{\cal F}_{0}\!\times\! {\cal F}_{1}\!\rightarrow\! X_{0}\!\times\! X_{1}~\ \ 
x\! =\!\varphi_{\varepsilon}(f_{\varepsilon})\ \ \hbox{\rm and}\ \ 
\forall t\!\in\! {\cal T}~\varphi (t)\!\in\!\Phi (t)\ \}.$$
As $\tilde\varphi (t_{0})\! =\! [\tilde\varphi_{0}(f_{0}),\tilde
\varphi_{1}(f_{1})]\!\in\!\Phi (t_{0})\cap (U_{0}\!\times\! U_{1})$ we get 
$\emptyset\!\not=\!\Phi (t_{0})\cap (U_{0}\!\times\! U_{1})\!\subseteq\! 
\overline{\Psi (t_{0})}^{\tau_{1}}\cap (U_{0}\!\times\! U_{1})$. By the 
separation theorem this implies that $\Psi (t_{0})\cap (U_{0}\!\times\! U_{1})$ 
is not empty and contains some point $(x_{0},x_{1})$. Fix $\varepsilon\!\in\! 2$. 
As $x_{\varepsilon}\!\in\! U_{\varepsilon}$ there is $\psi^\varepsilon\! :\! 
{\cal F}_{0}\!\times\! {\cal F}_{1}\!\rightarrow\! X_{0}\!\times\! X_{1}$ such 
that $x_{\varepsilon}\! =\!\psi^\varepsilon_{\varepsilon}(f_{\varepsilon})$ and 
$\psi^\varepsilon (t)\!\in\!\Phi (t)$, for each $t\!\in\! {\cal T}$.\bigskip

\noindent $\bullet$ If $e_{0}\!\not=\! e'_{0}\!\in\! {\cal F}_{0}$ and 
$[(\tilde e_{i},j_{i})]_{i\leq l}$ is a path in $G_{\cal T}$ with 
$(\tilde e_{0},j_{0})\! =\!\overline{e_{0}}$ and $(\tilde e_{l},j_{l})\! =\!
\overline{e'_{0}}$, then it is unique by Theorem I.2.5 in [B]. We call it 
$p_{e_{0},e'_{0}}$. We will define a partition of 
${\cal F}_{0}\!\times\! {\cal F}_{1}$. We put
$$\begin{array}{ll}
N\!\!\!\! & :=\{\ (e_{0},e_{1})\!\in\! {\cal F}_{0}\!\times\! {\cal F}_{1}\!
\setminus\{ t_{0}\}\mid (e_{0},e_{1})\!\notin\! {\cal T}\ \ \hbox{\rm or}\ \ 
[e_{0}\!\not=\! f_{0}\ \hbox{\rm and}\ p_{e_{0},f_{0}}\ \hbox{\rm does\ not\ 
exist}]\ \}\hbox{\rm ,}\cr & \cr
H\!\!\!\! & :=\{\ (e_{0},e_{1})\!\in\! {\cal T}\!\setminus \{ t_{0}\}\mid 
e_{0}\!\not=\! f_{0}\ \hbox{\rm and}\ p_{e_{0},f_{0}}(|p_{e_{0},f_{0}}|\! -\! 2)
\! =\!\overline{f_{1}}\ \}\hbox{\rm ,}\cr & \cr 
V\!\!\!\! & :=\{\ (e_{0},e_{1})\!\in\! {\cal T}\!\setminus \{ t_{0}\}\mid 
e_{0}\! =\! f_{0}\ \ \hbox{\rm or}\ \ [e_{0}\!\not=\! f_{0}\ \hbox{\rm and}
\ p_{e_{0},f_{0}}(|p_{e_{0},f_{0}}|\! -\! 2)\!\not=\!\overline{f_{1}}]\ \}.
\end{array}$$
The definition of $H$ means that if we view the graph $G_{\cal T}$ as $\cal T$ 
itself in the product ${\cal F}_{0}\!\times\! {\cal F}_{1}$ instead of seeing it 
in the sum ${\cal F}_{0}\!\oplus\! {\cal F}_{1}$, then the last edge in the path 
from $(e_{0},e_{1})$ to $t_{0}$ is horizontal (and vertical in $V$). So we 
defined a partition $(\{ t_{0}\},N,H,V)$ of ${\cal F}_{0}\!\times\! {\cal F}_{1}$.
\bigskip 

\noindent $\bullet$ Let us show that $\Pi_{{\cal F}_{\varepsilon}}[H]\cap
\Pi_{{\cal F}_{\varepsilon}}[V]\! =\!\emptyset$, for each $\varepsilon\!\in\! 2$.
\bigskip

 We may assume that $\varepsilon\! =\! 1$. We argue by contradiction. This gives 
$e_{1}\!\in\!\Pi_{{\cal F}_{1}}[H]\cap\Pi_{{\cal F}_{1}}[V]$, and also $e_{0}$ 
(resp., $e'_{0}$) such that $(e_{0},e_{1})\!\in\! H$ (resp., 
$(e'_{0},e_{1})\!\in\! V$). Note that $e_{0}\!\not=\! f_{0}$, and also that 
$e_{1}\!\not=\! f_{1}$ (by contradiction, we get $e'_{0}\!\not=\! f_{0}$ since 
$(e'_{0},e_{1})\!\not=\! t_{0}$, and 
$p_{e'_{0},f_{0}}\! =\! (\overline{e'_{0}},\overline{f_{1}},\overline{f_{0}})$, 
which is absurd). If $e'_{0}\! =\! f_{0}$, then 
$\overline{e_{1}}^\frown {p_{e_{0},f_{0}}}^\frown\overline{e_{1}}$ gives a cycle, 
which is absurd. If $e'_{0}\!\not=\! f_{0}$, then 
$\overline{e_{1}}^\frown p_{e_{0},f_{0}}$ and 
$\overline{e_{1}}^\frown p_{e'_{0},f_{0}}$ give two different pathes from 
$\overline{e_{1}}$ to $\overline{f_{0}}$, which is also absurd.\bigskip

\noindent $\bullet$ Now we can define 
$\psi_{\varepsilon}\! :\! {\cal F}_{\varepsilon}\!\rightarrow\! X_{\varepsilon}$. 
We put 
$$\begin{array}{ll}
\psi_{0}(e_{0})\!\!\!\! & :=\left\{\!\!\!\!\!\!
\begin{array}{ll}
& x_{0}\ \ \hbox{\rm if}\ \ e_{0}\! =\! f_{0}\hbox{\rm ,}\cr & \cr
& \psi^{1}_{0}(e_{0})\ \ \hbox{\rm if}\ \ e_{0}\!\in\!\Pi_{{\cal F}_{0}}[H]
\hbox{\rm ,}\cr & \cr
& \psi^{0}_{0}(e_{0})\ \ \hbox{\rm otherwise,}
\end{array}\right.\cr & \cr 
\psi_{1}(e_{1})\!\!\!\! & :=\left\{\!\!\!\!\!\!
\begin{array}{ll}
& x_{1}\ \ \hbox{\rm if}\ \ e_{1}\! =\! f_{1}\hbox{\rm ,}\cr & \cr
& \psi^{1}_{1}(e_{1})\ \ \hbox{\rm if}\ \ e_{1}\!\in\!\Pi_{{\cal F}_{1}}[H]\!
\setminus\!\{ f_{1}\}\hbox{\rm ,}\cr & \cr
& \psi^{0}_{1}(e_{1})\ \ \hbox{\rm otherwise.}
\end{array}\right.
\end{array}$$
Then we set $\psi (e_{0},e_{1})\! :=\! [\psi_{0}(e_{0}),\psi_{1}(e_{1})]$.

\vfill\eject

\noindent $\bullet$ It remains to see that $\psi (t)\!\in\!\Psi (t)$, for each 
$t\!\in\! {\cal T}$. Notice first that 
$\psi (t_{0})\! =\! (x_{0},x_{1})\!\in\!\Psi (t_{0})$. If 
$t\! :=\! (e_{0},e_{1})\!\in\! V$ and $e_{0}\!\not=\! f_{0}$, then we get
$$\psi (t)\! =\! [\psi_{0}(e_{0}),\psi_{1}(e_{1})]\! =\! 
[\psi^{0}_{0}(e_{0}),\psi^{0}_{1}(e_{1})]\! =\!\psi^{0}(t)\!\in\! 
\Phi (t)\! =\!\Psi (t).$$
Now if $t\!\in\! V$ and $e_{0}\! =\! f_{0}$, then we get
$$\psi (t)\! =\! [x_{0},\psi^{0}_{1}(e_{1})]\! =\! 
[\psi^{0}_{0}(f_{0}),\psi^{0}_{1}(e_{1})]\! =\!
[\psi^{0}_{0}(e_{0}),\psi^{0}_{1}(e_{1})]\! =\!
\psi^{0}(t)\!\in\!\Phi (t)\! =\!\Psi (t).$$
We argue similarly if $t\!\in\! H$.\bigskip

 If $t\!\in\! N\cap {\cal T}$, then $e_{0}\!\not=\! f_{0}$. If moreover 
$e_{1}\!\notin\! (\{ f_{1}\}\cup\Pi_{{\cal F}_{1}}[H])$, then we get
$$\psi (t)\! =\! [\psi_{0}(e_{0}),\psi_{1}(e_{1})]\! =\! 
[\psi^{0}_{0}(e_{0}),\psi^{0}_{1}(e_{1})]\! =\!\psi^{0}(t)\!\in\! 
\Phi (t)\! =\!\Psi (t).$$
If $e_{1}\! =\! f_{1}$, then 
$p_{e_{0},f_{0}}\! =\! (\overline{e_{0}},\overline{e_{1}},\overline{f_{0}})$ 
exists, which is absurd. If 
$e_{1}\!\in\!\Pi_{{\cal F}_{1}}[H]\!\setminus\{ f_{1}\}$, let 
$e'_{0}\!\in\! {\cal F}_{0}$ with $(e'_{0},e_{1})\!\in\! H$. The sequence 
$(\overline{e_{0}},\overline{e_{1}},\overline{e'_{0}},\ldots ,\overline{f_{1}},
\overline{f_{0}})$ shows that $p_{e_{0},f_{0}}$ exists, which is absurd again.
\bigskip

\noindent (b) Write ${\cal T}\!:=\!\{ t_{1},\ldots ,t_{n}\}$, and set 
$\Phi_{0}\!:=\!\overline{\Psi}$. We define ${\Phi_{j+1}\! :\! {\cal F}_{0}\!
\times\! {\cal F}_{1}\!\rightarrow\!\Ana (X_{0}\!\times\! X_{1})}$ as follows. We 
put $\Phi_{j+1}(t)\! :=\!\Phi_{j}(t)$ if $t\!\not=\! t_{j+1}$, and   
${\Phi_{j+1}(t_{j+1})\! :=\!\Psi (t_{j+1})}$, for $j\! <\! n$. The result now 
follows from an iterative application of (a).\hfill{$\square$}\bigskip

\noindent\bf Proof of Theorem 1.13.\rm ~We argue by contradiction. This gives a 
Borel set $S'$. Consider first that $A\! :=\! S'$ and $B\! :=\!\neg S'$. Then (b) 
holds with $u\! =\! v\! =\!\hbox{\rm Id}_{2^\omega}$. So (a) does not hold and 
$S'$ is not $\hbox{\rm pot}(\bormpx )$.\bigskip 

 Consider now that $A\! :=\! S$ and $B\! :=\!\lceil T\rceil\!\setminus\! S$, 
where $T$ and $S$ are given by Theorem 1.11. As (a) does not hold, (b) holds. 
This gives continuous maps $u$, $v$ with 
$$\begin{array}{ll} 
{\ \ \ }S'\!\!\!\! & \!\subseteq\! (u\!\times\! v)^{-1}(S)\!
\subseteq\! (u\!\times\! v)^{-1}(\lceil T\rceil )\hbox{\rm ,}\cr & \cr
\neg S'\!\!\!\! & \!\subseteq\! 
(u\!\times\! v)^{-1}(\lceil T\rceil\!\setminus\! S)\!\subseteq\! 
(u\!\times\! v)^{-1}(\lceil T\rceil ).
\end{array}$$
\bf Claim.\rm\ There is a Borel subset $A$ of $2^\omega$ with 
$S'\! =\! A\!\times\! 2^\omega$ or 
$S'\! =\! 2^\omega\!\times\! A$.\bigskip

\noindent $\bullet$ We argue by contradiction to prove the claim. There are  
$\alpha\!\in\! 2^\omega$, and $\beta\!\not=\!\beta'\!\in\! 2^\omega$ such that 
$(\alpha ,\beta )\!\in\! S'$ and 
$(\alpha ,\beta' )\!\notin\! S'$ (otherwise 
$A\! :=\! (S')^{0^\infty}\!\in\!\borel (2^\omega )$ and satisfies 
$S'\! =\! A\!\times\! 2^\omega$). Then 
$\big( u(\alpha ),v(\beta )\big)\!\in\! S$ and 
$\big( u(\alpha ),v(\beta' )\big)\!\notin\! S$, thus 
$v(\beta )\!\not=\! v(\beta')$.\bigskip 

\noindent $\bullet$ Note that $(\alpha',\beta )\!\in\! S'$, for each 
$\alpha'\!\in\! 2^\omega$. Indeed, we argue by contradiction. This 
gives $\alpha'$ with $\big( u(\alpha' ),v(\beta )\big)\!\notin\! S$. Thus 
$u(\alpha )\!\not=\! u(\alpha' )$, and $\big( u(\alpha ),v(\beta )\big)$, 
$\big( u(\alpha' ),v(\beta )\big)$, $\big( u(\alpha ),v(\beta' )\big)$, 
$\big( u(\alpha' ),v(\beta' )\big)$ are in $\lceil T\rceil$. Let 
$p\!\in\!\omega$ with 
$e_{0}\! :=\! u(\alpha )\lceil p\not= e'_{0}\! :=\! u(\alpha' )\lceil p\ 
\ \hbox{\rm and}\ \ e_{1}\! :=\! v(\beta )\lceil p\not= e'_{1}\! :=\! v(\beta' )
\lceil p.$ Then $(e_{0},e_{1})$, $(e'_{0},e_{1})$, $(e_{0},e'_{1})$, 
$(e'_{0},e'_{1})\!\in\! {\cal T}_{p}$, and the sequence $(\overline{e_{0}},
\overline{e_{1}},\overline{e'_{0}},\overline{e'_{1}},\overline{e_{0}})$ is a 
cycle, which is absurd.\bigskip

\noindent $\bullet$ Let $\gamma\!\in\! S'_{\alpha}$. We have 
$(\alpha',\gamma )\!\in\! S'$, for each $\alpha'\!\in\! 2^\omega$, as 
before. Conversely, assume that $(\alpha',\gamma )\!\in\! S'$. Then 
$\gamma\!\in\! S'_{\alpha}$, as before. Thus  
$S'\! =\! 2^\omega\!\times\! S'_{\alpha}$, which is absurd. 
This proves the claim.\hfill{$\diamond$}\bigskip

 Now the claim contradicts the fact that $S'$ is not 
$\hbox{\rm pot}(\bormpx )$.\hfill{$\square$}

\vfill\eject

\subsection{$\!\!\!\!\!\!$ The topologies.}\indent

 In this subsection we prove another result that will be used to show Theorem 
1.10. Some topologies are involved, and this is the place where we use some 
results in [Lo2].\bigskip

\noindent\bf Notation.\rm ~Let $X$, $Y$ be recursively presented Polish spaces.
\bigskip

\noindent $\bullet$ Recall the existence of $\Ca$ sets 
$W^{X}\!\subseteq\!\omega$, $C^{X}\!\subseteq\!\omega\!\times\! X$ with 
${\Borel (X)\! =\!\{ C^{X}_{n}\mid n\!\in\! W^{X}\}}$ and 
${\{(n,x)\!\in\!\omega\!\times\!\! X\mid n\!\in\! W^{X}\hbox{\rm and}~x\!\notin\! 
C_{n}^{X}\}\!\in\!\Ca (\omega\!\times\! X)}$ (see Theorem 3.3.1 in [H-K-Lo]).
\bigskip

\noindent $\bullet$ Set $\hbox{\rm pot}(\bormz)\! :=\!\borel (X)\!\times\!
\borel (Y)$ and, for $\xi\! <\!\omega^{\hbox{\rm CK}}_{1}$,
$$W^{X\times Y}_{\xi}\! :=\!\{ p\!\in\! W^{X\times Y}\mid 
C^{X\times Y}_{p}\!\in\!\hbox{\rm pot}(\bormxi)\}.$$ 
We also set ${W^{X\times Y}_{<\xi}\! :=\!\bigcup_{\eta <\xi}~
W^{X\times Y}_{\eta}}$.\bigskip

 The following result is essentially proved in [Lo2]. However, the statement is 
not in it, so we give a proof, which uses several statements in [Lo2]. 
Recall that $\tau_{1}$ is defined before Lemma 2.1.2.

\begin{them} (Louveau) Let $\xi\! <\!\omega^{\hbox{\it CK}}_{1}$, 
$X$, $Y$ be recursively presented Polish spaces. Then 
$W^{X\times Y}_{\xi}$ and $W^{X\times Y}_{<\xi}$ are $\Ca$. If moreover $A$, $B$ 
are disjoint $\Ana$ subsets of $X\!\times\! Y$, then the following are equivalent:

\noindent (a) The set $A$ is separable from $B$ by a $\hbox{\it pot}(\bormpx )$ 
set.

\noindent (b) The set $A$ is separable from $B$ by a 
$\Borel\cap\hbox{\it pot}(\bormpx )$ set.

\noindent (c) The set $A$ is separable from $B$ by a $\bormpx (\tau_{1})$ set.
\end{them}

\noindent\bf Proof.\rm ~By the second paragraph page 44 in [Lo2], $\borel (X)$ 
and $\borel (Y)$ are regular families (see Definition 2.7 in [Lo2] for the 
definition of a regular family). By Theorem 2.12 in [Lo2], the family 
$\Phi\! :=\!\hbox{\rm pot}(\bormz )$ is regular too. We define a sequence 
$(\Phi_{\xi})_{\xi <\omega_{1}^{\rm CK}}$ of families as follows (see Corollary 
2.10.(v) in [Lo2]):
$$\begin{array}{ll}\Phi_{0}\!\!\!\!\! 
& :=\!\Phi\hbox{\rm ,}\cr\Phi_{\xi +1}\!\!\!\!\! 
& :=\! (\Phi_{\xi})_{\sigma\hbox{\rm c}}\hbox{\rm ,}\cr\Phi_{\lambda}\!\!\!\!\! 
& :=\! \bigcup_{\xi <\lambda}\ \Phi_{\xi}\ \ \hbox{\rm if}\ \ 0\! 
<\!\lambda\! <\!\omega_{1}^{\hbox{\rm CK}}\ \hbox{\rm is\ a\ limit\ ordinal.}
\end{array}$$
By Corollary 2.10.(v) in [Lo2], $\Phi_{\xi}$ is a regular family for each 
$\xi\! <\!\omega_{1}^{\hbox{\rm CK}}$. In particular, the set $W_{\Phi_{\xi}}\! 
:=\!\{ p\!\in\! W^{X\times Y}\mid C_{p}^{X\times Y}\!\in\!\Phi_{\xi}\}$ is 
$\Ca (\omega )$. By Theorem 2.8 in [Lo2], the family $\Phi_{\xi +1}$ is a 
separating family (see Definition 2.1 in [Lo2] for the definition of a separating 
family), for each $\xi\! <\!\omega_{1}^{\hbox{\rm CK}}$. An easy induction on 
$\xi$ shows the following facts:
$$\begin{array}{ll}\Phi_{\xi}\!\!\!\!\! 
& =\!\hbox{\rm pot}(\bormxi )\ \ \hbox{\rm if}\ \ \xi\! <\!\omega\hbox{\rm ,}\cr
\Phi_{\xi}\!\!\!\!\! 
& =\!\bigcup_{\eta <\xi}\ \hbox{\rm pot}(\borme )\ \ \hbox{\rm if}\ \ 
0\! <\!\xi\! <\!\omega_{1}^{\hbox{\rm CK}}\ \hbox{\rm is\ a\ limit\ ordinal,}\cr
\Phi_{\xi +1}\!\!\!\!\! 
& =\!\hbox{\rm pot}(\bormxi )\ \ \hbox{\rm if}\ \ \omega\!\leq\!\xi\! <\!
\omega_{1}^{\hbox{\rm CK}}. 
\end{array}$$
This shows that $W_{\xi}^{X\times Y}\! =\! W_{\Phi_{\xi}}$ is $\Ca$ if 
$\xi\! <\!\omega$, $W_{\xi}^{X\times Y}\! =\! W_{\Phi_{\xi+1}}$ is $\Ca$ if 
$\omega\!\leq\!\xi\! <\!\omega_{1}^{\hbox{\rm CK}}$. If 
$0\! <\!\xi\! <\!\omega_{1}^{\hbox{\rm CK}}$ is a limit ordinal, then 
$W_{<\xi}^{X\times Y}\! =\! W_{\Phi_{\xi}}$ is $\Ca$.

\vfill\eject

\noindent (b) $\Rightarrow$ (c) follows from Theorem 3.4 in [Lo2].\bigskip

\noindent (c) $\Rightarrow$ (a) follows from the fact that ${\it\Delta}_{X}$ and 
${\it\Delta}_{Y}$ are Polish.\bigskip

\noindent (a) $\Rightarrow$ (b) Assume first that $\xi\! <\!\omega$. Then 
$\hbox{\rm pot}(\bormpx )\! =\! \Phi_{1+\xi}\! =\!\Phi_{\xi+1}$ is a separating 
family. So $A$ and $B$ are separable by a 
$\Borel\cap\Phi_{\xi+1}\! =\!\Borel\cap\hbox{\rm pot}(\bormpx )$ set. If 
$\omega\!\leq\!\xi\! <\!\omega_{1}^{\hbox{\rm CK}}$, then we use the fact that 
$\hbox{\rm pot}(\bormpx )\! =\!\hbox{\rm pot}(\bormxi )\! =\!\Phi_{\xi+1}$.
\hfill{$\square$}\bigskip

\noindent\bf Notation.\rm ~Let $X$, $Y$ be recursively presented Polish spaces.
\bigskip

\noindent $\bullet$ We will use the Gandy-Harrington topology 
${\it\Sigma}_{X}$ on $X$ generated by $\Ana (X)$. Recall that the set 
$\Omega_{X}\! :=\!\{x\!\in\! X\mid\omega_1^x\! =\!\omega^{\hbox{\rm CK}}_{1}\}$ is 
Borel and $\Ana$, that $[\Omega_{X},{\it\Sigma}_{X}]$ is a $0$-dimensional 
Polish space (the intersection of $\Omega_{X}$ with any nonempty $\Ana$ 
set is a nonempty clopen subset of $[\Omega_{X},{\it\Sigma}_{X}]$) (see [L1]).
\bigskip 

\noindent $\bullet$ Let $2\!\leq\!\xi\! <\!\omega^{\hbox{\rm CK}}_{1}$. The 
topology $\tau_{\xi}$ is generated by ${\Ana (X\!\times\! Y)\cap\bormlxi 
(\tau_{1})}$. We have the inclusion 
${\boraone(\tau_{\xi})\!\subseteq\!\boraxi (\tau_{1})}$, so that 
${\bormone(\tau_{\xi})\!\subseteq\!\bormxi (\tau_{1})}$. These topologies 
are similar to the ones considered in [Lo1] (see Definition 1.5).

\begin{lemm} Let $X$, $Y$ be recursively presented Polish 
spaces, and $\xi\! <\!\omega^{\hbox{\it CK}}_{1}$.\smallskip

\noindent (a) Fix $S\!\in\!\Ana (X\!\times\! Y)$. Then 
$\overline{S}^{\tau_{1+\xi}}\!\in\!\Ana (X\!\times\! Y)$.\smallskip

\noindent (b) Let $n\!\geq\! 1$, $1\!\leq\!\xi_{1}\! <\!\xi_{2}\! <\!\ldots\! <\!
\xi_{n}\!\leq\! 1\! +\!\xi$, and $S_{1}$, $\ldots$, $S_{n}$ be $\Ana$ sets. 
Assume that $S_{n'}\!\subseteq\!\overline{S_{n'+1}}^{\tau_{\xi_{n'}+1}}$ for 
$1\!\leq\! n'\! <\! n$. Then 
${S_{n}\cap\bigcap_{1\leq i<n}~\overline{S_{i}}^{\tau_{\xi_{i}}}}$ is 
$\tau_{1}$-dense in $\overline{S_{1}}^{\tau_{1}}$.\end{lemm}

\noindent\bf Proof.\rm ~(a) This is essentially proved in [Lo2] (see the proof of 
Theorem 2.8 in [Lo2]). We emphasize the fact that the analogous version of (a) in 
[Lo2] and the assertions of Theorem 2.2.1 are proved simultaneously by induction 
on $\xi$, and interact. Assume first that $\xi\! =\! 0$. Then 
$$\begin{array}{ll}
(x,y)\!\notin\!\overline{S}^{\tau_{1}}\!\!\!\! & \Leftrightarrow 
~\exists U\!\in\!\Borel (X)\ \exists V\!\in\!\Borel (Y)~~
(x,y)\!\in\! U\!\times\! V~\hbox{\rm and}~(U\!\times\! V)\cap S\! =\!\emptyset\cr   
& \Leftrightarrow  
~\exists m\!\in\! W^{X}\ \exists n\!\in\! W^{Y}\ \big(~C^{X}_{m}(x)
~~\hbox{\rm and}~~C^{Y}_{n}(y)~~\hbox{\rm and}~~\forall (x',y')\!\in\! 
X\!\times\! Y\cr
& \ \ \ \ \ \ [\big( m\!\in\! W^{X}~\hbox{\rm and}~x'\!\notin\! C_{m}^{X}\big)
~\hbox{\rm or}~\big( n\!\in\! W^{Y}~\hbox{\rm and}~y'\!\notin\! C_{n}^{Y}\big)
~\hbox{\rm or}~(x',y')\!\notin\! S]~\big).
\end{array}$$
So $\overline{S}^{\tau_{1}}\!\in\!\Ana (X\!\times\! Y)$. Now assume that 
$\xi\!\geq\! 1$. We have, by Theorem 2.2.1:
$$\begin{array}{ll}
(x,y)\!\notin\!\overline{S}^{\tau_{1+\xi}}\!\!\!\! & \Leftrightarrow 
~\exists T\!\in\!\Ana (X\!\times\! Y)\cap\bormlpxi (\tau_{1})~~
(x,y)\!\in\! T~\hbox{\rm and}~T\cap S\! =\!\emptyset\cr   
& \Leftrightarrow  
~\exists E\!\in\!\Borel (X\!\times\! Y)\cap \hbox{\rm pot}(\bormlpxi )~~
(x,y)\!\in\! E~\hbox{\rm and}~E\cap S\! =\!\emptyset\cr   
& \Leftrightarrow
~\exists m\!\in\! W^{X\times Y}_{<1+\xi}~~\big(~C^{X\times 
Y}_{m}(x,y)~~\hbox{\rm and}~~\forall (x',y')\!\in\! X\!\times\! Y\cr
& \ \ \ \ \ \ [\big( m\!\in\! W^{X\times Y}~\hbox{\rm and}~(x',y')\!\notin\! 
C_{m}^{X\times Y}\big)~\hbox{\rm or}~(x',y')\!\notin\! S]~\big).
\end{array}$$
By Theorem 2.2.1, $W^{X\times Y}_{<1+\xi}\!\in\!\Ca$ and we are done.\bigskip

\noindent (b) Let $U$ (resp., $V$) a $\Borel (X)$ (resp., $\Borel (Y)$) set with 
$\overline{S_{1}}^{\tau_{1}}\cap (U\!\times\! V)\!\not=\!\emptyset$. Then 
$S_{1}\cap (U\!\times\! V)\!\not=\!\emptyset$, which proves the desired property 
for $n\! =\! 1$. Then we argue inductively on $n$. So assume that the property is 
proved for $n$. We have $S_{n}\!\subseteq\!\overline{S_{n+1}}^{\tau_{\xi_{n}+1}}$, 
and $S_{n}\cap\bigcap_{1\leq i<n}~\overline{S_{i}}^{\tau_{\xi_{i}}}\cap 
(U\!\times\! V)\!\not=\!\emptyset$, by induction assumption. Thus 
$\overline{S_{n+1}}^{\tau_{\xi_{n}+1}}\cap\bigcap_{1\leq i\leq n}~
\overline{S_{i}}^{\tau_{\xi_{i}}}\cap (U\!\times\! V)\!\not=\!\emptyset$. As 
$\bigcap_{1\leq i\leq n}~\overline{S_{i}}^{\tau_{\xi_{i}}}\cap (U\!\times\! V)$ 
is $\boraone (\tau_{\xi_{n}+1})$, we get $S_{n+1}\cap\bigcap_{1\leq i\leq n}~
\overline{S_{i}}^{\tau_{\xi_{i}}}\cap (U\!\times\! V)\!\not=\!\emptyset$.
\hfill{$\square$}

\vfill\eject

\subsection{$\!\!\!\!\!\!$ Representation of Borel sets.}\indent

 Now we come to the representation theorem of Borel sets by G. 
Debs and J. Saint Raymond (see [D-SR]). It specifies the classical 
result of Lusin asserting that any Borel set in a Polish space is the 
bijective continuous image of a closed subset of the Baire space. The following 
definitions can be found in [D-SR]:

\begin{defin} (Debs-Saint Raymond) Let $a$ be a finite set. A partial order 
relation $R$ on $a^{<\omega}$ is a $\underline{tree~relation}$ if, for 
${t\!\in\! a^{<\omega}}$,\smallskip

\noindent (a) $\emptyset ~R~t$.\smallskip

\noindent (b) The set $P_{R}(t)\! :=\!\{s\!\in\! a^{<\omega}\mid s~R~t\}$ is finite 
and linearly ordered by $R$.\smallskip

 For instance, the non strict extension relation $\prec$ is a tree relation.
\smallskip

\noindent $\bullet$ Let $R$ be a tree relation. An $\underline{R\! -\! branch}$ 
is an $\subseteq$-maximal subset of $a^{<\omega}$ linearly ordered by $R$. We 
denote by $[R]$ the set of all infinite $R$-branches.\smallskip

 We equip $(a^{<\omega})^\omega$ with the product of the discrete topology on 
$a^{<\omega}$. If $R$ is a tree relation, the space 
$[R]\!\subseteq\! (a^{<\omega})^\omega$ is equipped with the topology induced by 
that of $(a^{<\omega})^\omega$. The map 
$\theta\! :\! a^\omega\!\rightarrow\! [\prec ]$ defined by 
$\theta (\gamma )\! :=\! [\gamma\lceil j]_{j\in\omega}$ is an homeomorphism.\smallskip

\noindent $\bullet$ Let $R$, $S$ be tree relations with $R\!\subseteq\! S$. The 
$\underline{canonical\ map}$ $\Pi\! :\! [R]\!\rightarrow\! [S]$ is defined by
$$\Pi (A)\! :=\!\hbox{\it the unique $S$-branch containing $A$.}$$
$\bullet$ Let $S$ be a tree relation. We say that $R\!\subseteq\! S$ is 
$\underline{distinguished}$ in $S$ if
$$\forall s,t,u\!\in\! a^{<\omega}\ \ \left.
\begin{array}{ll}
& s~S~t~S~u\cr & \cr
& \ \ s~R~u
\end{array}\right\} ~\Rightarrow ~s~R~t.$$
For example, let $C$ be a closed subset of $a^\omega$, and define: 
$$s~R~t\ \Leftrightarrow\ s\!\prec\! t\ \hbox{\rm and}\ N_{t}\cap 
C\!\not=\!\emptyset .$$ 
Then $R$ is distinguished in $\prec$. In this case, the distinction expresses 
the fact that ``when we leave the closed set, it is for ever".\smallskip

\noindent $\bullet$ Let $\eta\! <\!\omega_{1}$. A family 
$(R^{(\rho )})_{\rho\leq\eta}$ of tree relations is a 
$\underline{resolution~family}$ if:\smallskip

\noindent (a) $R^{(\rho +1)}$ is a distinguished subtree of $R^{(\rho )}$, for 
all $\rho\! <\!\eta$.\smallskip

\noindent (b) $R^{(\lambda )}\! =\!\bigcap_{\rho <\lambda}~R^{(\rho )}$, for 
all limit $\lambda\!\leq\!\eta$.\end{defin}

 We will use the following extension of the property of distinction:

\begin{lemm} Let $\eta\! <\!\omega_{1}$, 
$(R^{(\rho )})_{\rho\leq\eta}$ a resolution family with $R^{(0)}=\ \prec$, and 
$\rho\! <\!\eta$. Assume that $s\!\prec\! s'~R^{(\rho )}~s''$ and 
$s~R^{(\rho +1)}~s''$. Then $s~R^{(\rho +1)}~s'$.\end{lemm}

\noindent\bf Proof.\rm ~We argue by induction on $\rho$. Assume 
that the property is proved for $\mu\! <\!\rho$. As $s'~R^{(\rho)}~s''$ 
and $R^{(\rho +1)}$ is distinguished in $R^{(\rho)}$ we have $s~R^{(\rho 
+1)}~s'$.\hfill{$\square$}

\vfill\eject

\noindent\bf Notation.\rm ~Let $\eta\! <\!\omega_{1}$, 
$(R^{(\rho )})_{\rho\leq\eta}$ a resolution family with $R^{(0)}=\ \prec$, 
$\rho\!\leq\!\eta$ and $z\!\in\! a^{<\omega}\!\setminus\!\{\emptyset\}$. We set
$$z^{\rho}\! :=\! z~\lceil ~\hbox{\rm max}\{ r\! <\! |z|\mid z\lceil r~R^{(\rho )}~z
\}.$$
We enumerate $\{ z^{\rho}\mid\rho\!\leq\!\eta\}$ by 
$\{ z^{\xi_{i}}\mid 1\!\leq\! i\!\leq\! n\}$, where $1\!\leq\! n\!\in\!\omega$ and 
${\xi_{1}\! <\!\ldots\! <\!\xi_{n}\! =\!\eta}$. We can write 
$z^{\xi_{n}}\!\prec_{\not=}\! 
z^{\xi_{n-1}}\!\prec_{\not=}\!\ldots\!\prec_{\not=}\! z^{\xi_{2}}\!\prec_{\not=}
\! z^{\xi_{1}}\!\prec_{\not=}\! z$. By Lemma 2.3.2 we have 
$z^{\xi_{i+1}}~R^{(\xi_{i}+1)}~z^{\xi_{i}}$ for each $1\!\leq\! i\! <\! n$.

\begin{lemm} Let $\eta\! <\!\omega_{1}$, 
$(R^{(\rho )})_{\rho\leq\eta}$ a resolution family with $R^{(0)}=\ \prec$, 
$z\!\in\! a^{<\omega}\!\setminus\!\{\emptyset\}$ and $1\!\leq\! i\! <\! n$.
\smallskip

\noindent (a) Set $\eta_{i}\! :=\!\{\rho\!\leq\!\eta\mid z^{\xi_{i}}\!\prec\! 
z^{\rho}\}$. Then $\eta_{i}$ is a successor ordinal.\smallskip

\noindent (b) We may assume that 
$z^{\xi_{i}+1}\!\prec_{\not=}\! z^{\xi_{i}}$.\end{lemm}

\noindent\bf Proof.\rm ~(a) First notice that $\eta_{i}$ is an ordinal. Note that 
$\xi_{i}\! +\! 1\!\leq\!\eta_{i}\!\leq\!\eta\! +\! 1$. We argue by contradiction, 
so that $\eta_{i}\!\leq\!\eta$. Let $\xi_{i}\!\leq\!\rho\! <\!\eta_{i}$. Then we 
have $z^{\xi_{i}}\! =\! z^{\rho}$, $z^{\xi_{i}}~R^{(\rho )}~z$, 
$z^{\xi_{i}}~R^{(\eta_{i})}~z$, and $z^{\xi_{i}}\!\prec\! z^{\eta_{i}}$. 
As $\eta_{i}\!\leq\!\eta$, we get $\eta_{i}\!\in\!\eta_{i}$, which is 
absurd.\bigskip

\noindent (b)  So we can write $\eta_{i}\! =\!\nu_{i}\! +\! 1$. Note that 
$z^{\nu_{i}}\! =\! z^{\xi_{i}}$ since $\xi_{i}\!\leq\!\nu_{i}$. If 
$\nu_{i}\! +\! 1\!\leq\!\eta$ 
we get $z^{\nu_{i}+1}\!\prec_{\not=}\! z^{\nu_{i}}$, so we may assume that 
$\xi_{i}\! =\!\nu_{i}$. If $\nu_{i}\! +\! 1\! =\!\eta\! +\! 1$ we get 
$\nu_{i}\! =\!\eta$ and 
$z^{\xi_{i}}\! =\! z^{\nu_{i}}\! =\! z^{\eta}\! =\! z^{\xi_{n}}$, which is 
absurd.\hfill{$\square$}\bigskip

 The following is part of Theorem I-6.6 in [D-SR].

\begin{them} (Debs-Saint Raymond) Let $\eta\! <\!\omega_{1}$, $E$ be 
a $\bormep$ subset of $[\prec ]$. Then there is a resolution family 
$(R^{(\rho )})_{\rho\leq\eta}$ with:\smallskip

\noindent (a) $R^{(0)}=\ \prec$.\smallskip

\noindent (b) The canonical map $\Pi\! :\! [R^{(\eta )}]\!\rightarrow\! [\prec ]$ 
is a bijection.\smallskip

\noindent (c) The set $\Pi^{-1}(E)$ is a closed subset of $[R^{(\eta )}]$.
\end{them}

 Now we come to the actual proof of Theorem 1.10. 

\subsection{$\!\!\!\!\!\!$ Proof of Theorem 1.10.}\indent

\begin{them} Let $T$ be a tree with acyclic levels, 
$\xi\! <\!\omega^{\hbox{\it CK}}_{1}$ such that $1\! +\!\xi$ is a successor 
ordinal, $S\!\in\!\borapx (\lceil T\rceil)$, $X$, $Y$ recursively 
presented Polish spaces, and $A$, $B$ disjoint $\Ana$ subsets of $X\!\times\! Y$. 
Then one of the following holds:\smallskip

\noindent (a) $\overline{A}^{\tau_{1+\xi}}\cap B\! =\!\emptyset$.\smallskip

\noindent (b) There are $u\! :\! 2^\omega\!\rightarrow\! X$ and 
$v\! :\! 2^\omega\!\rightarrow\! Y$ continuous with  
$S\!\subseteq\! (u\!\times\! v)^{-1}(A)$ and 
$\lceil T\rceil\!\setminus\! S\!\subseteq\! (u\!\times\! v)^{-1}(B)$.
\end{them}

\noindent\bf Proof.\rm ~Fix $\eta\! <\!\omega^{\hbox{\rm CK}}_{1}$ with 
$1\! +\!\xi\! =\!\eta\! +\! 1$.\bigskip

\noindent $\bullet$ We identify $(2\!\times\! 2)^{Q}$ with $2^{Q}\!\times\! 2^{Q}$, for 
$Q\!\leq\!\omega$. With the notation of Definition 2.3.1 and $a\! :=\! 2\!\times\! 2$, we get 
$E\!:=\!\theta [\lceil T\rceil\!\setminus\! S]\!\in\!\bormep ([\prec ])$. 
Theorem 2.3.4 provides a resolution family. We put 
$$D\! :=\!\{(s,t)\!\in\! T\mid\exists\gamma\!\in\!\Pi^{-1}(E)\ \ (s,t)\!\in\!\gamma\}.$$ 
For example, we may assume that ${(\emptyset ,\emptyset )\!\in\! D}$.

\vfill\eject

\noindent $\bullet$ We set ${N\! :=\!\overline{A}^{\tau_{1+\xi}}\cap B}$. 
Applying Lemma 2.2.2.(a), we see that $N$ is $\Ana$. We assume that $N$ is not 
empty. Recall that $[\Omega_{X\times Y},{\it\Sigma}_{X\times Y}]$ is a Polish 
space (see the notation before Lemma 2.2.2). We fix a complete metric $d$ (resp., 
metrics $\delta_{X}$, $\delta_{Y}$) on 
$[\Omega_{X\times Y},{\it\Sigma}_{X\times Y}]$ (resp., $X$, $Y$ equipped with the 
initial topologies).\bigskip

\noindent $\bullet$ We construct $(x_{s})_{s\in\Pi_{0}[T]}\!\subseteq\! X$, 
$(y_{t})_{t\in\Pi_{1}[T]}\!\subseteq\! Y$,     
${(U_{(s,t)})_{(s,t)\in T}\!\subseteq\!\Ana (X\!\times\! Y)}$ with:
$$
\begin{array}{ll} 
& (\hbox{\rm i})~~~~(x_{s},y_{t})\!\in\! U_{(s,t)}\!\subseteq\!
\Omega_{X\times Y}.\cr & \cr
& (\hbox{\rm ii})~~~\hbox{\rm diam}_{d}(U_{(s,t)})\!\leq\! 2^{-|s|},~ 
\delta_{X}(x_{s},x_{s\varepsilon })\!\leq\! 2^{-|s|},~ 
\delta_{Y}(y_{t},y_{t\varepsilon })\!\leq\! 2^{-|t|}.\cr & \cr
& (\hbox{\rm iii})~~U_{(s,t)}\!\subseteq\! N~\hbox{\rm if}~(s,t)\!\in\! D.\cr & \cr
& (\hbox{\rm iv})~~~U_{(s,t)}\!\subseteq\! A~\hbox{\rm if}~(s,t)\!\notin\! D.\cr & \cr
& (\hbox{\rm v})~~~~[1\!\leq\!\rho\!\leq\!\eta ~\ \hbox{\rm and}~\ (s,t)~R^{(\rho 
)}~(s',t')]~\Rightarrow ~U_{(s',t')}\!\subseteq\!
\overline{U_{(s,t)}}^{\tau_{\rho }}.\cr & \cr
& (\hbox{\rm vi})~~~[\big( (s,t)\!\in\! D~\Leftrightarrow ~(s',t')\!\in\! D\big)~~
\hbox{\rm and}~~(s,t)~R^{(\eta )}~(s',t')]~\Rightarrow ~U_{(s',t')}\!\subseteq\! 
U_{(s,t)}.
\end{array}$$
$\bullet$ Let us show that this construction is sufficient to get the 
theorem. If $(\alpha ,\beta )\!\in\!\lceil T\rceil$, then we can define  
${(j_{i})_{i\in\omega}\! :=\! (j^{\alpha ,\beta}_{i})_{i\in\omega}}$ by 
${\Pi^{-1}\big( [(\alpha ,\beta )\lceil j]_{j\in\omega}\big)\! =\! 
[(\alpha ,\beta )\lceil j_{i}]_{i\in\omega}}$, where $j_{i}\! <\! j_{i+1}$. In 
particular, we have 
${(\alpha ,\beta )\lceil j_{i}~R^{(\eta )}~(\alpha ,\beta )\lceil 
j_{i+1}}$. We have the following:
$$\begin{array}{ll} 
(\alpha ,\beta )\!\in\! S\!\!\! 
& \Leftrightarrow ~\theta (\alpha ,\beta )\! =\! 
[(\alpha ,\beta )\lceil j]_{j\in\omega}\!\notin\! E\ 
\Leftrightarrow ~[(\alpha ,\beta )\lceil j_{i}]_{i\in\omega}
\!\notin\!\Pi^{-1}(E)\cr & \cr
& \Leftrightarrow ~\exists i_{0}\!\in\!\omega~\forall i\!\geq\! i_{0}~~
(\alpha ,\beta )\lceil j_{i}\!\notin\! D
\end{array}$$
since $\Pi^{-1}(E)$ is a closed subset of $[R^{(\eta )}]$. Similarly, 
$(\alpha ,\beta )\!\in\!\lceil T\rceil\!\setminus\! S$ is equivalent to the 
existence of $i_{0}\!\in\!\omega$ such that $(\alpha ,\beta )\lceil j_{i}\!\in\!  D$ for each 
$i\!\geq\! i_{0}$ (with $i_{0}\! =\! 0$).\bigskip

 Therefore ${U_{(\alpha ,\beta )\lceil j_{i+1}}\!\subseteq\! U_{(\alpha ,\beta 
)\lceil j_{i}}\!\subseteq\!\Omega_{X\times Y}}$ if $i\!\geq\! i_{0}$ and 
$(\alpha ,\beta )\!\in\!\lceil T\rceil$. Thus 
$(U_{(\alpha,\beta )\lceil j_{i}})_{i\geq i_{0}}$ is a decreasing sequence of 
nonempty clopen subsets of $[\Omega_{X\times Y},d]$ whose diameters 
tend to $0$. Therefore 
${\{F(\alpha ,\beta )\}\! =\!\bigcap_{i\geq i_{0}}~U_{(\alpha,\beta 
)\lceil j_{i}}}$ 
defines $F(\alpha ,\beta )$ in $\Omega_{X\times Y}$. Note that 
$F(\alpha ,\beta )$ is the limit of the sequence 
$\big( (x_{\alpha\lceil j_{i}},y_{\beta\lceil j_{i}})\big)_{i\in\omega}$.\bigskip

 Let $\alpha\!\in\!\Pi_{0}(\lceil T\rceil)$, and $\beta_{\alpha}$ such 
that $(\alpha ,\beta_{\alpha})\!\in\!\lceil T\rceil$. We set 
$u(\alpha )\! :=\!\Pi_{X}\big( F(\alpha ,\beta_{\alpha})\big)$. Note 
that $u(\alpha )$ is the limit of some subsequence of $(x_{\alpha\lceil 
i})_{i\in\omega}$, by continuity of the projection. As 
$\delta_{X}(x_{s},x_{s\varepsilon})\!\leq\! 2^{-|s|}$, $u(\alpha )$ is 
also the limit of $(x_{\alpha\lceil i})_{i\in\omega}$. Thus $u(\alpha )$ does 
not depend on the choice of $\beta_{\alpha}$. This also shows that $u$ is 
continuous on 
$\Pi_{0}(\lceil T\rceil)$. As $\Pi_{0}(\lceil T\rceil)$ is a closed 
subset of $2^\omega$, we can find a continuous retraction $r_{0}$ from 
$2^\omega$ onto $\Pi_{0}(\lceil T\rceil)$ (see Proposition 2.8 in [K]). We set 
$u(\alpha )\! :=\! u\big( r_{0}(\alpha )\big)$, so that $u$ is continuous 
on $2^\omega$.\bigskip

 Similarly, we define a continuous map $v\! :\! 2^\omega\!\rightarrow\! Y$ such 
that $v(\beta )$ is the limit of $(y_{\beta\lceil i})_{i\in\omega}$ if $\beta$ is in 
$\Pi_{1}(\lceil T\rceil)$. This implies that 
$F(\alpha ,\beta )\! =\!\big( u(\alpha ),v(\beta )\big)$ if 
$(\alpha ,\beta )\!\in\!\lceil T\rceil$.\bigskip

 If $(\alpha ,\beta )\!\in\! S$ (resp., $\lceil T\rceil\!\setminus\! S$), then 
$F(\alpha ,\beta )\!\in\! A$ (resp., $N$). This shows that 
${S\!\subseteq\! (u\!\times\! v)^{-1}(A)}$ and 
$\lceil T\rceil\!\setminus\! S\!\subseteq\! (u\!\times\! v)^{-1}(B)$.

\vfill\eject
 
\noindent $\bullet$ So let us show that the construction is possible. 
Fix $(x_{\emptyset},y_{\emptyset})\!\in\! N\cap\Omega_{X\times Y}$, which is 
not empty since $N\!\not=\!\emptyset$ is $\Ana$. Then we choose 
$U_{(\emptyset ,\emptyset )}\!\in\!\Ana$ with diameter at most $1$ with 
${(x_{\emptyset},y_{\emptyset})\!\in\! U_{(\emptyset ,\emptyset 
)}\!\subseteq\! N\cap\Omega_{X\times Y}}$. Assume that 
$(x_{s})_{|s|\leq p}$, $(y_{t})_{|t| \leq p}$, $(U_{(s,t)})_{|s|\leq p}$ 
satisfying conditions (i)-(vi) have been constructed, which is the case for 
$p\! =\! 0$.\bigskip

\noindent - Let $s\!\in\!\Pi_{0}[T]\cap 2^p$ (resp., 
$t\!\in\!\Pi_{1}[T]\cap 2^p$), and $X_{s}$ (resp., $Y_{t}$) be a 
$\Borel$ neighborhood of $x_{s}$ (resp., $y_{t}$) with $\delta_{X}$-diameter 
(resp., $\delta_{Y}$-diameter) at most $2^{-p}$.\bigskip 

\noindent - If $we\! :=\! (s\varepsilon ,t\varepsilon ')\!\in\! T\cap 
(2\!\times\! 2)^{p+1}$ $(w\! :=\! (s,t)\!\in\! (2\!\times\! 2)^{p}$ and 
$e\! :=\! (\varepsilon ,\varepsilon ')\!\in\! 2\!\times\! 2)$, then we set
$$(we)^{\eta +1}\! :=\!\left\{\!\!\!\!\!\!
\begin{array}{ll} 
& (we)^{\eta}\ \ \hbox{\rm if\ there\ is}\ \ r\!\leq\! p\ \ \hbox{\rm with}\ \ 
[~w\lceil r\!\in\! D~\Leftrightarrow ~we\!\in\! D~]\ \ \hbox{\rm and}\ \ 
w\lceil r~R^{(\eta )}~we\hbox{\rm ,}\cr 
& we\ \ \hbox{\rm otherwise.}
\end{array}\right.$$
Note that $(we)^{\eta}\!\in\! D$ if $we\!\in\! D$, so that $we\!\notin\! D$ if 
$(we)^{\eta +1}\! =\! we$. Note also the equivalence between the fact that 
$we\!\in\! D$, and the fact that $(we)^{\eta +1}\!\in\! D$. Indeed, we may assume 
that $we\!\notin\! D$ and $(we)^{\eta +1}\! =\! (we)^{\eta}$. So that there is 
$r\!\leq\! p$ with $w\lceil r\!\notin\! D$ and $w\lceil r~R^{(\eta )}~we$. By 
Lemma 2.3.2 we have $w\lceil r~R^{(\eta )}~(we)^{\eta}$, so that 
$(we)^{\eta +1}\! =\! (we)^{\eta}\!\notin\! D$. The conclusions in the assertions 
(a) and (b) in the following claim do not really depend on their respective 
assumptions, but we will use these assertions later in this form.\bigskip

\noindent\bf Claim.\rm ~Assume that $\eta\! >\! 0$.

\noindent (a) $A\cap\bigcap_{1\leq\rho\leq\eta}~
\overline{U_{(we)^{\rho}}}^{\tau_{\rho}}\cap (X_{s}\!\times\! Y_{t})$ is 
$\tau_{1}$-dense in $\overline{U_{(we)^1}}^{\tau_{1}}\cap 
(X_{s}\!\times\! Y_{t})$ if $(we)^{\eta +1}\! =\! we$.

\noindent (b) $U_{(we)^{\eta}}\cap\bigcap_{1\leq\rho <\eta}~
\overline{U_{(we)^{\rho}}}^{\tau_{\rho}}\cap (X_{s}\!\times\! Y_{t})$ is 
$\tau_{1}$-dense in $\overline{U_{(we)^1}}^{\tau_{1}}\cap 
(X_{s}\!\times\! Y_{t})$ if ${(we)^{\eta +1}\!\not=\! we}$.\bigskip  

 Indeed, we use the notation before Lemma 2.3.3 with $z\! :=\! we$. By Lemma 
2.3.3 we may assume that $z^{\xi_{i}+1}\!\prec_{\not=}\! z^{\xi_{i}}$ if 
$1\!\leq\! i\! <\! n$. We set ${S_{i}\! :=\! U_{z^{\xi_{i}}}}$, for 
${1\!\leq\!\xi_{i}\!\leq\!\eta}$. We have 
$S_i\!\subseteq\!\overline{S_{i+1}}^{\tau_{\xi_i+1}}$, for 
$1\!\leq\!\xi_{i}\! <\!\eta$, by induction assumption, since 
$z^{\xi_{i+1}}\ R^{(\xi_i+1)}\ z^{\xi_i}$. Moreover, the inclusion  
$S_n\!\subseteq\!\overline{A}^{\tau_{\eta +1}}$ holds. Thus 
$A\cap\bigcap_{1\leq\xi_i\leq\eta}~\overline{U_{(we)^{\xi_i}}}^{\tau_{\xi_i}}\cap 
(X_{s}\!\times\! Y_{t})$ 
(respectively, ${U_{(we)^{\eta}}\cap\bigcap_{1\leq\xi_i 
<\eta}~\overline{U_{(we)^{\xi_i}}}^{\tau_{\xi_i}}\cap 
(X_{s}\!\times\! Y_{t})}$) is $\tau_{1}$-dense in the set 
$\overline{U_{(we)^1}}^{\tau_{1}}\cap (X_{s}\!\times\! Y_{t})$ if 
${(we)^{\eta +1}\! =\! we}$ (respectively, $(we)^{\eta +1}\!\not=\! we$), 
by Lemma 2.2.2.(b). But if $1\!\leq\!\rho\!\leq\!\eta$, then there is $1\!\leq\! 
i\!\leq\! n$ with $(we)^{\rho}\! =\! (we)^{\xi_i}$. And $\rho\!\leq\!\xi_i$ since 
we have $(we)^{\xi_{i}+1}\!\prec_{\not=}\! (we)^{\xi_{i}}$ if ${1\!\leq\! i\! <\! n}$. 
Thus we are done since 
$\bigcap_{1\leq\rho\leq\eta}~\overline{U_{(we)^{\rho}}}^{\tau_{\rho}}\! =\! 
\bigcap_{1\leq\xi_i\leq\eta}~\overline{U_{(we)^{\xi_i}}}^{\tau_{\xi_i}}$ 
and $U_{(we)^{\eta}}\cap\bigcap_{1\leq\rho 
<\eta}~\overline{U_{(we)^{\rho}}}^{\tau_{\rho}}\! =\! 
U_{(we)^{\eta}}\cap\bigcap_{1\leq\xi_i 
<\eta}~\overline{U_{(we)^{\xi_i}}}^{\tau_{\xi_i}}$.\hfill{$\diamond$}\bigskip 

\noindent - Let ${\cal F}_{0}\!:=\! {\cal F}_{1}\!:=\! 2^{p+1}$, 
${\cal T}\!:=\! T\cap ({\cal F}_{0}\!\times\! {\cal F}_{1})$, $\Psi\! :\! {\cal 
F}_{0}\!\times\! {\cal F}_{1}\!\rightarrow\!\Ana (X\!\times\! Y)$ defined on 
${\cal T}$ by
$$\Psi (we)\! :=\!\left\{\!\!\!\!\!\!
\begin{array}{ll} 
& A\cap\bigcap_{1\leq\rho\leq\eta}~\overline{U_{(we)^{\rho}}}^{\tau_{\rho}}\cap
(X_{s}\!\times\! Y_{t})\cap\Omega_{X\times Y}~\hbox{\rm if}~(we)^{\eta +1}\! =\! 
we,\cr 
& U_{(we)^{\eta}}\cap\bigcap_{1\leq\rho <\eta}~\overline{U_{(we)^{\rho}}}^
{\tau_{\rho}}\cap (X_{s}\!\times\! Y_{t})~\hbox{\rm if}~(we)^{\eta +1}\!\not=\! 
we.
\end{array}\right.$$
By the claim, $\Psi (we)$ is $\tau_{1}$-dense in 
$\overline{U_{(we)^1}}^{\tau_{1}}\cap (X_{s}\!\times\! Y_{t})$ if $\eta\! 
>\! 0$. As $(we)^1\!\prec\! w\!\prec\! we$ and $R^{(1)}$ is distinguished in 
$\prec$ we get $(we)^1\ R^{(1)}\ w$ and 
$U_w\!\subseteq\!\overline{U_{(we)^1}}^{\tau_1}$, by induction 
assumption. Thus $\overline{U_{w}}^{\tau_1}\cap (X_s\!\times\! Y_t)\!\subseteq\!
\overline{U_{(we)^1}}^{\tau_1}\cap (X_s\!\times\! Y_t)\!\subseteq\!
\overline{\Psi}(we)$. Thus $(x_{s},y_{t})$ is in 
$U_{w}\cap (X_{s}\!\times\! Y_{t})\!\subseteq\!\overline{\Psi}(we)$ (even if 
$\eta\! =\! 0$). Therefore $\overline{\Psi}$ admits a $\pi$-selector on ${\cal T}$. 
By Lemma 2.1.2, $\Psi$ admits a $\pi$-selector $\psi$ on ${\cal T}$. We set 
$x_{s\varepsilon}\! :=\!\psi_{0}(s\varepsilon )$, 
$y_{t\varepsilon'}\! :=\!\psi_{1}(t\varepsilon')$, and choose $\Ana$ sets 
$U_{we}$ with $d$-diameter at most $2^{-p-1}$ such that 
$\psi (we)\!\in\! U_{we}\!\subseteq\!\Psi (we)$. This finishes the proof 
since $(s,t)\ R^{(\rho )}\ we$ and 
$(s,t)\!\not=\! we$ imply that $(s,t)\ R^{(\rho )}\ (we)^\rho\ R^{(\rho 
)}\ we$, by Lemma 2.3.2.\hfill{$\square$}

\vfill\eject

 Now we come to the limit case. We need some more definitions that can be found 
in [D-SR].

\begin{defin} (Debs-Saint Raymond) Let $a$ be a finite set.\smallskip

\noindent $\bullet$ Let $R$ be a tree relation on $a^{<\omega}$. If 
$t\!\in\! a^{<\omega}$, then $h_{R}(t)$ is the number of strict $R$-predecessors 
of $t$. So we have ${h_{R}(t)\! =\!\hbox{\it Card}\big( P_R(t)\big)\! -\! 1}$.
\smallskip

\noindent $\bullet$ Let $\xi\! <\!\omega_{1}$ be an infinite limit ordinal. We 
say that a resolution family $(R^{(\rho)})_{\rho\leq\xi}$ is 
$\underline{uniform}$ if
$$\forall k\!\in\!\omega\ \exists\eta_k\! <\!\xi\ 
\forall s,t\!\in\! a^{<\omega}\ \ 
[\hbox{\it min}\big( h_{R^{(\xi)}}(s),h_{R^{(\xi)}}(t)\big)\!\leq\! k\ 
\hbox{\it and}\ s\ R^{(\eta_k)}\ t]\ \Rightarrow\ s\ R^{(\xi)}\ t.$$
We may (and will) assume that $\eta_k\!\geq\! 1$.\end{defin}

 The following is part of Theorem I-6.6 in [D-SR].
 
\begin{them} (Debs-Saint Raymond) Let $\xi\! <\!\omega_{1}$ be an 
infinite limit ordinal, $E$ a $\bormxi$ subset of $[\prec ]$. Then there is a 
uniform resolution family $(R^{(\rho )})_{\rho\leq\xi}$ with:\smallskip

\noindent (a) $R^{(0)}=\ \prec$.\smallskip

\noindent (b) The canonical map $\Pi\! :\! [R^{(\xi )}]\!\rightarrow\! [\prec ]$ 
is a bijection.\smallskip

\noindent (c) The set $\Pi^{-1}(E)$ is a closed subset of $[R^{(\xi )}]$.
\end{them}

\begin{them} (Debs-Lecomte) Let $T$ be a tree with acyclic levels, 
$\xi\! <\!\omega^{\hbox{\it CK}}_{1}$ an infinite limit ordinal, 
$S\!\in\!\boraxi (\lceil T\rceil)$, $X$, $Y$ recursively presented Polish 
spaces, and $A$, $B$ disjoint $\Ana$ subsets of $X\!\times\! Y$. Then one of the 
following holds:\smallskip

\noindent (a) $\overline{A}^{\tau_{\xi}}\cap B\! =\!\emptyset$.\smallskip

\noindent (b) There are $u\! :\! 2^\omega\!\rightarrow\! X$ and 
$v\! :\! 2^\omega\!\rightarrow\! Y$ continuous with 
$S\!\subseteq\! (u\!\times\! v)^{-1}(A)$ and 
$\lceil T\rceil\!\setminus\! S\!\subseteq\! (u\!\times\! v)^{-1}(B)$.
\end{them}

\noindent\bf Proof.\rm ~Let us indicate the differences with the 
proof of Theorem 2.4.1.\bigskip

\noindent $\bullet$ The set ${E\!:=\!\theta [\lceil T\rceil\!\setminus\! S]}$ is 
$\bormxi ([\prec ])$. Theorem 2.4.3 provides a uniform resolution family.
\bigskip

\noindent $\bullet$ If $w\!\in\! (2\!\times\! 2)^{<\omega}$ then we set
$$\eta (w)\! :=\!\hbox{\rm max}\{\eta_{h_{R^{(\xi )}}(w')+1}\mid w'\!\prec\! w\}.$$
Note that $\eta (w')\!\leq\!\eta (w)$ if $w'\!\prec\! w$.\bigskip

\noindent $\bullet$ Conditions (v) and (vi) become\bigskip

\leftline{$
\begin{array}{ll}
& (\hbox{\rm v}')~~~~[1\!\leq\!\rho\!\leq\!\eta (s,t)~\ \hbox{\rm and}~\ 
(s,t)~R^{(\rho)}~(s',t')]~\Rightarrow 
~U_{(s',t')}\!\subseteq\!\overline{U_{(s,t)}}^{\tau_{\rho }}.\cr & \cr
& (\hbox{\rm vi}')~~~[\big( (s,t)\!\in\! D~\Leftrightarrow ~(s',t')\!\in\! D\big)
~~\hbox{\rm and}~~(s,t)~R^{(\xi )}~(s',t')]~\Rightarrow ~U_{(s',t')}\!\subseteq\! 
U_{(s,t)}.
\end{array}$}\bigskip

\noindent $\bullet$ If $we\! :=\! (s\varepsilon ,t\varepsilon ')\!\in\! T\cap 
(2\!\times\! 2)^{p+1}$, then we set
$$(we)^{\xi +1}\! :=\!\left\{\!\!\!\!\!\!
\begin{array}{ll} 
& (we)^{\xi}\ \ \hbox{\rm if\ there\ is}\ \ r\!\leq\! p\ \ \hbox{\rm with}\ \ 
[~w\lceil r\!\in\! D~\Leftrightarrow ~we\!\in\! D~]\ \ \hbox{\rm and}\ \ 
w\lceil r~R^{(\xi )}~we\hbox{\rm ,}\cr & \cr 
& we\ \ \hbox{\rm otherwise.}
\end{array}\right.$$
Note that $we\!\notin\! D$ if $(we)^{\xi +1}\! =\! we$. Note also the 
equivalence between the fact that $we\!\in\! D$ and 
the fact that $(we)^{\xi +1}\!\in\! D$.

\vfill\eject

\noindent\bf Claim 1.\rm ~Assume that $(we)^\rho\!\not=\! (we)^\xi$. Then 
$\rho\! +\! 1\!\leq\!\eta\big( (we)^{\rho +1}\big)$.\bigskip 

 We argue by contradiction. We get
$$\rho\! +\! 1\! >\!\rho\!\geq\!\eta\big( (we)^{\rho +1}\big)\!\geq\!
\eta_{h_{R^{(\xi )}}( (we)^{\xi})+1}\! =\!\eta_{h_{R^{(\xi )}}(we)}.$$
As $(we)^\rho\ R^{(\rho)}\ we$ we get $(we)^\rho\ R^{(\xi)}\ we$ and 
$(we)^\rho\! =\! (we)^\xi$, which is absurd.\hfill{$\diamond$}\bigskip

 Note that $\xi_{n-1}\! <\!\xi_{n-1}\! +\! 1\!\leq\!\eta\big( 
(we)^{\xi_{n-1}+1}\big)\!\leq\!\eta (we)$. Thus 
$(we)^{\eta (we)}\! =\! (we)^{\xi}$.\bigskip

\noindent\bf Claim 2.\rm ~(a) $A\cap\bigcap_{1\leq\rho\leq\eta (we)}~
\overline{U_{(we)^{\rho}}}^{\tau_{\rho}}\cap (X_{s}\!\times\! Y_{t})$ is 
$\tau_{1}$-dense in $\overline{U_{(we)^1}}^{\tau_{1}}\cap 
(X_{s}\!\times\! Y_{t})$ if $(we)^{\xi +1}\! =\! we$.

\noindent (b) $U_{(we)^{\xi}}\cap\bigcap_{1\leq\rho <\eta (we)}~
\overline{U_{(we)^{\rho}}}^{\tau_{\rho}}\cap (X_{s}\!\times\! Y_{t})$ is 
$\tau_{1}$-dense in $\overline{U_{(we)^1}}^{\tau_{1}}\cap 
(X_{s}\!\times\! Y_{t})$ if ${(we)^{\xi +1}\!\not=\! we}$.\bigskip 

 Indeed, we set $S_{i}\! :=\! U_{z^{\xi_{i}}}$, for $1\!\leq\!\xi_{i}\!\leq\!\xi$. 
By Claim 1 we can apply Lemma 2.2.2.(b) and we are done.\hfill{$\diamond$}
\bigskip

\noindent $\bullet$ Let ${\cal F}_{0}\!:=\! {\cal F}_{1}\!:=\! 2^{p+1}$, 
${\cal T}\!:=\! T\cap ({\cal F}_{0}\!\times\! {\cal F}_{1})$, $\Psi\! :\! 
{\cal F}_{0}\!\times\! {\cal F}_{1}\!\rightarrow\!\Ana (X\!\times\! Y)$ defined 
on ${\cal T}$ by
$$\Psi (we)\! :=\!\left\{\!\!\!\!\!\!
\begin{array}{ll} & 
A\cap\bigcap_{1\leq\rho\leq\eta(we)}~\overline{U_{(we)^{\rho}}}^{\tau_{\rho}}
\cap (X_{s}\!\times\! Y_{t})\cap\Omega_{X\times Y}\ \hbox{\rm if}~(we)^{\xi +1}\! 
=\! we,\cr & \cr
& U_{(we)^{\xi}}\cap\bigcap_{1\leq\rho <\eta (we)}~\overline{U_{(we)^{\rho}}}^
{\tau_{\rho}}\cap (X_{s}\!\times\! Y_{t})\ \ \hbox{\rm if}~(we)^{\xi +1}\!\not=\! 
we.
\end{array}\right.$$
We conclude as in the proof of Theorem 2.4.1, using the facts that 
$\eta_k\!\geq\! 1$ and $\eta (.)$ is increasing.\hfill{$\square$}\bigskip

\noindent\bf Proof of Theorem 1.10.\rm ~
We may assume that $\xi\! <\!\omega^{\hbox{\rm CK}}_{1}$, $X$, $Y$ are 
recursively presented, and $A,B$ are $\Ana$. We assume 
that $A$ is not separable from $B$ by a $\hbox{\rm pot}(\bormpx )$ set, and set 
${N\! :=\!\overline{A}^{\tau_{1+\xi}}\cap B}$. Then $N$ is not empty since 
${\bormone(\tau_{1+\xi})\!\subseteq\!\bormpx (\tau_{1})}\!\subseteq\!  
\hbox{\rm pot}(\bormpx )$. So (b) holds, by Theorems 2.4.1 and 2.4.3.\bigskip

 So (a) or (b) holds. If $D\!\in\!\hbox{\rm pot}(\bormpx )$ separates $A$ from $B$ and 
(b) holds, then $S\!\in\!\hbox{\rm pot}(\bormpx )$, since 
${S\! =\! (u\!\times\! v)^{-1}(D)\cap\lceil T\rceil}$,  which is absurd.\hfill{$\square$}

\section{$\!\!\!\!\!\!$ Proof of Theorem 1.11.}

 We have seen that we cannot have a reduction on the whole product in Theorem 
1.13. We have seen that it is possible to have it on the set of branches of some 
tree with acyclic levels. We now build an example of such a tree. This tree has 
to be small enough since we cannot have a reduction on the whole product. But as 
the same time it has to be big enough to ensure the existence of complicated sets, 
as in the statement of Theorem 1.11.\bigskip

\noindent\bf Notation.\rm ~Let 
$\varphi\! :\!\omega\!\rightarrow\!\omega^2$ be the natural bijection. 
More precisely, we set, for $q\!\in\!\omega$,
$$M(q)\! :=\!\hbox{\rm max}\{ m\!\in\!\omega\mid\Sigma_{k\leq m}~k\!\leq\! q\}.$$
Then we define 
$\varphi (q)\! =\!\big((q)_{0},(q)_{1}\big)\! :=\!\big(M(q)\! -\! q\! +\! 
(\Sigma_{k\leq M(q)}~k),q\! -\! (\Sigma_{k\leq M(q)}~k)\big)$. One can check that 
$<\! n,p\! >:=\!\varphi^{-1}(n,p)\! =\! 
(\Sigma_{k\leq n+p}~k)\! +\! p$. More concretely, we get
$$\varphi[\omega ]\! =\!\{(0,0),(1,0),(0,1),(2,0),(1,1),(0,2),\ldots\}.$$

\vfill\eject

\begin{defi} We say that 
$E\!\subseteq\!\bigcup_{q\in\omega}~2^q\!\times\! 2^q$ is a 
$\underline{ test}$ if:\medskip

\noindent (a) $\forall q\!\in\!\omega~\exists !(s_{q},t_{q})\!\in\! 
E\!\cap\! (2^q\!\times\! 2^q)$.\medskip

\noindent (b) $\forall m,p\!\in\!\omega~\forall u\!\in\! 
2^{<\omega}~\exists v\!\in\! 2^{<\omega}~(s_{p}0uv,t_{p}1uv)\!\in\! E$ and 
$(|t^{}_{p}1uv|\! -\! 1)_{0}\! =\! m$.\medskip

\noindent (c) $\forall n\! >\! 0~\exists q\! <\! n~\exists w\!\in\! 
2^{<\omega}~~s_{n}\! =\! s_{q}0w~\hbox{\rm and}~t_{n}\! =\! t_{q}1w$.\medskip

 We will call $T$ the tree generated by a test 
$E\! =\!\{(s_{q},t_{q})\mid q\!\in\!\omega\}$: 
$$T\! :=\!\{(s,t)\!\in\! 2^{<\omega}\!\times\! 2^{<\omega}\mid s\! =\! 
t\! =\!\emptyset ~~\hbox{\rm or}~~\exists q\!\in\!\omega ~\exists w\!\in\! 
2^{<\omega}~~(s,t)\! =\! (s_{q}0w,t_{q}1w)\}.$$\end{defi}

 The uniqueness condition in (a) and condition (c) ensure that $T$ is small 
enough, and also the acyclicity. The existence condition in (a) and 
condition (b) ensure that $T$ is big enough. More specifically, if 
$X$ is a Polish space and $\sigma$ a finer Polish topology on $X$, 
then there is a dense $G_{\delta}$ subset of $X$ on which the two 
topologies coincide. The first part of condition (b) ensures the 
possibility to get inside the square of a dense $G_{\delta}$ subset of 
$2^\omega$. The examples in Theorem 1.11 are build using the 
examples in [Lo-SR]. Conditions on verticals are involved, and the 
second part of condition (b) gives a control on the choice of verticals.

\begin{prop} The tree $T$ associated with a test is a tree with acyclic 
levels.\end{prop}

\noindent\bf Proof.\rm ~Fix $p\!\in\!\omega$. Let us show that $G_{{\cal T}_{p}}$ 
is acyclic. We argue by contradiction. Let $(\tilde e_i,j_{i})_{i\leq l}$ be a 
cycle in $G_{{\cal T}_{p}}$, and $n\! <\! p$ maximal such that the sequence 
$(\tilde e_i(n))_{i\leq l}$ is not constant. There is $i_1$ minimal with 
$\tilde e_{i_1}(n)\!\not=\!\tilde e_{i_1+1}(n)$. We have 
${\tilde e_{i_1}(n)\! =\!\tilde e_0(n)\! =\!\tilde e_l(n)}$. There is 
$i_2\! >\! i_1\! +\! 1$ minimal with 
${\tilde e_{i_1+1}(n)\!\not=\!\tilde e_{i_2}(n)}$. Then 
${\tilde e_{i_1}(n)\! =\!\tilde e_{i_2}(n)}$, and in fact 
$\tilde e_{i_1}\! =\!\tilde e_{i_2}$ because of the uniqueness condition in (a), 
and $\tilde e_{i_1+1}\! =\!\tilde e_{i_2-1}$. If $j_{i_{1}}\! =\! j_{i_{2}}$, 
then $i_1\! =\! 0$ and $i_2\! =\! l$. But 
$j_{i_{1}+1}\! =\! 1\! -\! j_{i_{1}}\! =\! 1\! -\! j_{i_{2}}\! =\! j_{i_{2}-1}$, 
which is absurd. If $j_{i_{1}}\!\not=\! j_{i_{2}}$, then for example 
$j_{i_{1}}\! =\! 0\! =\! 1\! -\! j_{i_{2}}$. If $p\! >\! 0$, then 
$\tilde e_{i_1}(0)\! =\! 0\! =\! 1\! -\!\tilde e_{i_2}(0)$, which contradicts 
$\tilde e_{i_1}\! =\!\tilde e_{i_2}$. If $p\! =\! 0$, then 
$\tilde e_{0}\! =\!\emptyset\! =\!\tilde e_{2}$, which is absurd.
\hfill{$\square$}\bigskip

\noindent\bf Notation.\rm ~Let 
$\psi\! :\!\omega\!\rightarrow\! 2^{<\omega}$ be the natural bijection 
($\psi (0)\! =\!\emptyset$, $\psi (1)\! =\! 0$, $\psi (2)\! =\! 1$, 
$\psi (3)\! =\! 0^2$, $\psi (4)\! =\! 01$, $\psi (5)\! =\! 10$, 
$\psi (6)\! =\! 1^2$, $\ldots$). Note that $|\psi (q)|\!\leq\! q$.

\begin{lem} There exists a test.\end{lem}

\noindent\bf Proof.\rm ~We set $s^{}_{0}\! =\! t^{}_{0}\! :=\!\emptyset$, and
$$\begin{array}{ll} 
s^{}_{q+1} 
& \!\!\!\!\! := s^{}_{[(q)_{1}]_{0}}~0~
\psi\big([(q)_{1}]_{1}\big)~0^{q-[(q)_{1}]_{0}-
|\psi([(q)_{1}]_{1})|},\cr & \cr t^{}_{q+1} 
& \!\!\!\!\! := t^{}_{[(q)_{1}]_{0}}~1~
\psi\big([(q)_{1}]_{1}\big)~0^{q-[(q)_{1}]_{0}-
|\psi([(q)_{1}]_{1})|}. 
\end{array}$$
Note that $(q)_{0}\! +\! (q)_{1}\! =\! M(q)\!\leq\!\Sigma_{k\leq M(q)}~k\!\leq\! 
q$, so that $s^{}_{q}, t^{}_{q}$ are well defined and we have the 
equality $|s^{}_{q}|\! =\! |t^{}_{q}|\! =\! q$, by 
induction on $q$. It remains to check that condition (b) in the 
definition of a test is fullfilled. Set $n\! :=\!\psi^{-1}(u)$, $r\! := 
<\! p,n\! >$ and $q\! :=<\! m,r\! >$. It remains to put $v\!:=\! 0^{q-p-|u|}$: 
${(s^{}_{p}0uv,t^{}_{p}1uv)\! =\! 
(s^{}_{q+1},t^{}_{q+1})}$.\hfill{$\square$}\bigskip

 Now we come to the lemma crucial for the proof of Theorem 1.11.

\vfill\eject

\noindent\bf Notation.\rm ~(a) We define 
$p\! :\!\omega^{<\omega}\!\setminus\!\{\emptyset\}\!\rightarrow\!\omega$. We 
define $p(s)$ by induction on $|s|$:
$$p(s)\! :=\!\left\{\!\!\!\!\!\!
\begin{array}{ll} 
& s(0)~\hbox{\rm if}~|s|\! =\! 1,\cr & \cr
& <\! p\big(s\lceil (|s|\! -\! 1)\big),s(|s|\! -\! 1)\! >~\hbox{\rm otherwise.}
\end{array}
\right.$$
Note that $p_{\vert\omega^n}\! :\!\omega^n\!\rightarrow\!\omega$ is a bijection, 
for each $n\!\geq\! 1$.\bigskip

\noindent (b) The map $\Delta\! :\! 2^{\omega}\!\times\! 
2^{\omega}\!\rightarrow\! 2^{\omega}$ is the symmetric difference. So, for 
$m\!\in\!\omega$, 
$$\Delta (\alpha ,\beta )(m)\! =\! (\alpha\Delta\beta )(m)\! =\! 1~
\Leftrightarrow ~\alpha (m)\!\not=\!\beta (m).$$

\begin{lem} Let $G$ be a dense $G_{\delta}$ subset of $2^\omega$, and $T$ the 
tree associated with a test. Then there are $\alpha_{0}\!\in\! G$ and 
$f\! :\! 2^\omega\!\rightarrow\! G$ continuous such that, for each 
$\alpha\!\in\! 2^\omega$,

\noindent (a) $\big(\alpha_{0},f(\alpha )\big)\!\in\!\lceil T\rceil$.

\noindent (b) For each $t\!\in\!\omega^{<\omega}$, and each 
$m\!\in\!\omega$,

(i) $\alpha\big( p(tm)\big)\! =\! 1~\Rightarrow ~\exists m'\!\in\!\omega ~~
\big(\alpha_{0}\Delta f(\alpha )\big)\big( 
p(tm')\! +\! 1\big)\! =\! 1$.

(ii) $\big(\alpha_{0}\Delta f(\alpha )\big)\big( 
p(tm)\! +\! 1\big)\! =\! 1~\Rightarrow ~
\exists m'\!\in\!\omega ~~\alpha\big( p(tm')\big)\! =\! 1$.\end{lem}

\noindent\bf Proof.\rm ~Let $(O_{q})$ be a sequence of dense 
open subsets of $2^\omega$ with $G\! =\!\bigcap_{q}~O_{q}$. By density we 
get: 
$\forall q,l\!\in\!\omega~\exists u_{q,l}\!\in\! 2^{<\omega}~\forall 
s\!\in\! 2^l~~N_{su_{q,l}}\!\subseteq\! O_{q}$.\bigskip

\noindent $\bullet$ We construct finite approximations of $\alpha_{0}$ and $f$. 
The idea is to linearize the binary tree $2^{<\omega}$. So we will 
use the bijection $\psi$ defined before Lemma 3.3. To construct 
$f(\alpha )$ we have to imagine, for each length $l$, the different possibilities 
for $\alpha\lceil l$. More precisely, we construct subsequences of 
$2^{<\omega}$, namely $(v_{w})_{w\in 2^{<\omega}}$, $(s_{w})_{w\in 2^{<\omega}}$ 
and $(t_{w})_{w\in 2^{<\omega}}$, satisfying the following conditions:
$$\begin{array}{ll} 
& (1)~\ (s_{w},t_{w})\!\in\! 
E\!\setminus\!\{ (\emptyset ,\emptyset)\},~\hbox{\rm and}~
(|t_{w}|\! -\! 1)_{0}\! =\! (|w|)_{0},~\hbox{\rm for~each}~w\!\in\! 2^{<\omega}.
\cr & \cr
& (2)~\ \left\{
\begin{array}{ll} 
s_{\emptyset} & \!\!\!\!\! =\! 0~u_{0,1}~v_{\emptyset},\cr & \cr 
s_{w\varepsilon } & \!\!\!\!\! =\! s_{\psi (\psi^{-1}(w\varepsilon )-1)}~0~
u_{\psi^{-1}(w\varepsilon ),|s_{\psi (\psi^{-1}(w\varepsilon 
)-1)}|+1}~v_{w\varepsilon }.
\end{array}\right.\cr & \cr 
& (3)~\ \left\{
\begin{array}{ll} 
t_{\emptyset} & \!\!\!\!\! =\! 1~u_{0,1}~v_{\emptyset},\cr & \cr 
t_{w\varepsilon } & \!\!\!\!\! =\! t_{w}\ \varepsilon\  
[{^\frown}_{~\psi^{-1}(w)<i<\psi^{-1}(w\varepsilon )}~u_{i,|s_{\psi 
(i-1)}|+1}~v_{\psi (i)}~0]\ u_{\psi^{-1}(w\varepsilon ),
|s_{\psi (\psi^{-1}(w\varepsilon )-1)}|+1}~v_{w\varepsilon }. 
\end{array}\right.
\end{array}$$
We show the existence of the three subsequences inductively on $\psi^{-1}(w)$. We choose 
$v_{\emptyset}\!\in\! 2^{<\omega}$ with 
$(0~u_{0,1}~v_{\emptyset},1~u_{0,1}~v_{\emptyset})\!\in\! E$ and 
$(|1~u_{0,1}~v_{\emptyset}|\! -\! 1)_{0}\! =\! 0$. Assume that 
$(v_{w})_{\psi^{-1}(w)\leq r}$, $(s_{w})_{\psi^{-1}(w)\leq r}$, 
$(t_{w})_{\psi^{-1}(w)\leq r}$ satisfying properties (1)-(3) have 
been constructed, which is the case for $r\! =\! 0$.\bigskip 

 Fix $w\!\in\! 2^{<\omega}$ and $\varepsilon\!\in\! 2$ with 
$\psi (r\! +\! 1)\! =\! w\varepsilon$. We choose 
$v_{w\varepsilon }\!\in\! 2^{<\omega}$ such that $(s_{w\varepsilon 
},t_{w\varepsilon })\!\in\! E$ and 
$(|t_{w\varepsilon }|\! -\! 1)_{0}\! =\! (|w|\! +\! 1)_{0}$. Let us show that 
this is possible. We want that\bigskip

\leftline{$(s_{\psi (\psi^{-1}(w\varepsilon 
)-1)}~0~u_{\psi^{-1}(w\varepsilon ),|s_{\psi (\psi^{-1}(w\varepsilon 
)-1)}|+1}~v_{w\varepsilon }~,~t_{w}~\varepsilon ~u_{\psi^{-1}(w)+1,|s_{w}|+1}~
v_{\psi (\psi^{-1}(w)+1)}\ 0\ldots$}\bigskip
\rightline{$u_{\psi^{-1}(w\varepsilon )-1,
|s_{\psi (\psi^{-1}(w\varepsilon )-2)}|+1}~
v_{\psi (\psi^{-1}(w\varepsilon )-1)}~0~u_{\psi^{-1}(w\varepsilon ),
|s_{\psi (\psi^{-1}(w\varepsilon )-1)}|+1}~v_{w\varepsilon })\!\in\! 
E.$}\bigskip\noindent
It is enough to see that $(s_{\psi (\psi^{-1}(w\varepsilon )-1)}~0,
t_{w}~\varepsilon \ldots v_{\psi (\psi^{-1}(w\varepsilon )-1)}~0)\!\in\! T$. 

\vfill\eject

 But
$$\begin{array}{ll} 
& s_{\psi (\psi^{-1}(w\varepsilon )-1)}~0\cr & \cr 
=\!\! & s_{\psi (\psi^{-1}(w\varepsilon )-2)}~0~u_{\psi^{-1}(w\varepsilon )-1,
|s_{\psi (\psi^{-1}(w\varepsilon )-2)}|+1}~
v_{\psi (\psi^{-1}(w\varepsilon )-1)}~0\cr  & \cr 
=\!\! & \ldots\cr  & \cr 
=\!\! & s_{w}~0~u_{\psi^{-1}(w)+1,|s_{w}|+1}~v_{\psi 
(\psi^{-1}(w)+1)}~0\ldots\ 
u_{\psi^{-1}(w\varepsilon )-1,|s_{\psi (\psi^{-1}(w\varepsilon 
)-2)}|+1}~v_{\psi (\psi^{-1}(w\varepsilon )-1)}~0.
\end{array}$$
We are done since $(s_{w},t_{w})\!\in\! E$.\bigskip

\noindent $\bullet$ So this defines sequences $(v_{w})_{w\in 2^{<\omega}}$, 
$(s_{w})_{w\in 2^{<\omega}}$ and $(t_{w})_{w\in 2^{<\omega}}$. As 
$s_{\psi (q)}\!\prec_{\not=}\! s_{\psi (q+1)}$ we can define 
${\alpha_{0}\! :=\!\hbox{\rm sup}_{q}~s_{\psi (q)}}$. Similarly, we set 
$f(\alpha )\! :=\!\hbox{\rm sup}_{m}~t_{\alpha\lceil m}$, and $f$ is 
continuous.\bigskip

\noindent $\bullet$ Let us show that $\alpha_{0}\!\in\! G$. By 
definition of $s_{w\varepsilon }$ we get 
${s_{\psi (q)}0u_{q+1,|s_{\psi (q)}|+1}\!\prec\! s_{\psi (q+1)}}$, 
for each $q$. This implies that 
$\alpha_{0}\!\in\!\bigcap_{q}~O_{q}\! =\! G$ since 
$0u_{0,1}\!\prec\!\alpha_{0}$.\bigskip

\noindent $\bullet$ Now fix $\alpha\!\in\! 2^\omega$. Let us show that 
$f(\alpha )\!\in\! G$. Fix $q\!\in\!\omega$, and $m\!\in\!\omega$ such that
$$\psi^{-1}(\alpha\lceil m)\! <\! q\! +\! 1\!\leq\!
\psi^{-1}\big(\alpha\lceil (m\! +\! 1)\big).$$
Again it is enough to show the existence of ${s\!\in\! 2^{<\omega}}$ with 
$su_{q+1,|s|}\!\prec\! t_{\alpha\lceil (m+1)}$. Set 
$$s\! :=\! t_{\alpha\lceil m}~\alpha (m)~u_{\psi^{-1}(\alpha\lceil m)+1,
|s_{\alpha\lceil m}|+1}~v_{\psi (\psi^{-1}(\alpha\lceil m)+1)}~0~\ldots ~
u_{q,|s_{\psi (q-1)}|+1}~v_{\psi (q)}~0.$$ 
By definition of $t_{\alpha\lceil (m+1)}$ we have 
$su_{q+1,|s_{\psi (q)}|+1}\!\prec\! t_{\alpha\lceil (m+1)}$. But the 
construction of $t_{w\varepsilon }$ shows that $|s_{\psi (q)}|\! +\! 1\! =\! |s|$. 
So $s$ is suitable.\bigskip

\noindent (a) Moreover, $\big(\alpha_{0}, f(\alpha 
)\big)\!\in\!\lceil T\rceil$. 
Indeed, fix $r\!\in\!\omega$. There is $m\!\in\!\omega$ with 
${l\! :=\! |t_{\alpha\lceil m}|\!\geq\! r}$. We get 
$\big(\alpha_{0},f(\alpha )\big)\lceil l\! =\!
\big(s_{\alpha\lceil m},t_{\alpha\lceil m}\big)\!\in\! E\!\subseteq\! T$. 
Thus $\big(\alpha_{0}, f(\alpha )\big)\lceil r\!\in\! T$, and 
$\big(\alpha_{0}, f(\alpha )\big)$ is in $\lceil T\rceil$.\bigskip

\noindent (b).(i) We set ${w\!:=\!\alpha\lceil p(tm)}$, so that 
${t_{w}1\!\prec\! t_{w1}\! =\! t_{\alpha\lceil [p(tm)+1]}\!\prec\! f(\alpha )}$. As 
${(|t_{w}|\! -\! 1)_{0}\! =\! p(t)}$, there is $m'$ with ${|t_{w}|\! =\! p(tm')\! +\! 1}$. But 
$s_{w}0\!\prec\! s_{\psi(\psi^{-1}(w)+1)}$, so that 
${\alpha_{0}(|t_{w}|)\!\not=\! f(\alpha )(|t_{w}|)}$.\bigskip

\noindent (b).(ii) First notice that the only coordinates where 
$\alpha_{0}$ and $f(\alpha )$ can differ are $0$ and the $|t_{\alpha\lceil q}|$'s. 
Therefore there is an integer $q$ with 
$p(tm)\! +\! 1\! =\! |t_{\alpha\lceil q}|$. In 
particular $(|t_{\alpha\lceil q}|\! -\! 1)_{0}\! =\! p(t)$ and 
$(q)_{0}\! =\! p(t)$. 
Thus there is $m'$ with $q\! =\! p(tm')$. We have 
$\alpha_{0}(|t_{\alpha\lceil q}|)\! =\! 0\!\not=\! f(\alpha )(|t_{\alpha\lceil 
q}|)\! =\!\alpha (q)$. So 
$\alpha\big(p(tm')\big)\! =\! 1$.\hfill{$\square$}\bigskip

 Now we come to the existence of complicated sets, as in the statement of Theorem 
1.11.\bigskip

\noindent\bf Notation.\rm ~In [Lo-SR], Lemma 3.3, the map 
$\rho_{0}\! :\! 2^\omega\!\rightarrow\! 2^\omega$ defined as follows is 
introduced:
$$\rho_{0}(\varepsilon )(i)\! :=\!\left\{\!\!\!\!\!\!
\begin{array}{ll}
& 1\ \ \hbox{\rm if}\ \ \varepsilon (<\! i,j\! >)\! =\! 0,\ \hbox{\rm for\ each\ }
j\!\in\!\omega ,\cr &Ê\cr
& 0\ \ \hbox{\rm otherwise.}
\end{array}\right.$$

\vfill\eject

 In this paper, $\rho_{0}^{\xi}\! :\! 2^\omega\!\rightarrow\! 2^\omega$ is also 
defined for $\xi\! <\!\omega_{1}$ as follows, by induction on $\xi$ (see the 
proof of Theorem 3.2). We put:

\noindent - ${\rho^0_{0}\! :=\!\hbox{\rm Id}_{2^\omega}}$.

\noindent - $\rho^{\eta +1}_{0}\!:=\!\rho^{}_{0}\circ\rho^{\eta}_{0}$.

\noindent - If $\lambda\! >\! 0$ is limit, fix 
$(\xi^\lambda_{k})\!\subseteq\!\lambda\!\setminus\!\{ 0\}$ with  
$\Sigma_{k}~\xi^\lambda_{k}\! =\!\lambda$. For 
$\varepsilon \!\in\! 2^\omega$ and $k\!\in\!\omega$ we define 
$(\varepsilon )^k\!\in\! 2^\omega$ 
by $(\varepsilon )^k(i)\! :=\!\varepsilon (i\! +\! k)$. We also define 
$\rho^{(k,k+1)}_{0}\! :\! 2^\omega\!\rightarrow\! 2^\omega$ by
$$\rho^{(k,k+1)}_{0}(\varepsilon )(i)\! :=\!\left\{\!\!\!\!\!\!
\begin{array}{ll} 
& \varepsilon (i)~\hbox{\rm if}~i\! <\! k,\cr &Ê\cr
& \rho_{0}^{\xi^\lambda_{k}}\big( (\varepsilon )^k\big)(i\! -\! k)~
\hbox{\rm if}~i\!\geq\! k.
\end{array}\right.$$ 
We set $\rho^{(0,k+1)}_{0}\! :=\!\rho^{(k,k+1)}_{0}\circ\rho^{(k-1,k)}_{0}\circ
\ldots\circ\rho^{(0,1)}_{0}$ and 
$\rho^\lambda_{0}(\varepsilon )(k)\! :=\!\rho^{(0,k+1)}_{0}(\varepsilon )(k)$.\bigskip 

 The set ${H_{1+\xi}\!:=\! (\rho^{\xi}_{0})^{-1}(\{ 0^\infty\})}$ is also 
introduced, and the authors show that $H_{1+\xi}$ is 
$\bormpx\!\setminus\!\borapx$ (see Theorem 3.2).\bigskip 

\noindent $\bullet$ The map ${\cal S}\! :\! 2^\omega\!\rightarrow\! 2^\omega$ is 
the shift map: ${\cal S}(\alpha )(m)\! :=\!\alpha (m\! +\! 1)$.\bigskip

\noindent $\bullet$ Let $T$ be the tree generated by a test. We put, 
for $\xi\! <\!\omega_{1}$,
$$S_{1+\xi}:=\{ (\alpha ,\beta )\!\in\! 2^\omega\!\times\! 
2^\omega\mid 
(\alpha ,\beta )\!\in\!\lceil T\rceil ~\hbox{\rm and}~{\cal S}(\alpha\Delta\beta )\!
\notin\! H_{1+\xi}\}.$$

\begin{thm} Let $\xi\! <\!\omega_{1}$. The set 
$\lceil T\rceil\!\setminus\! S_{1+\xi}$ is $\bormpx (2^\omega\!\times\! 
2^\omega)\!\setminus\hbox{\it pot}(\borapx )$, and 
$S_{1+\xi}$ is not $\hbox{\it pot}( \bormpx )$.\end{thm}

\noindent\bf Proof.\rm ~As $H_{1+\xi}$ is $\bormpx$ and $\Delta$, ${\cal S}$ 
are continuous, $\lceil T\rceil\!\setminus\! S_{1+\xi}$ is $\bormpx 
(2^\omega\!\times\! 2^\omega)$.\bigskip

\noindent $\bullet$ Let $G$ be a dense $G_{\delta}$ subset of 
$2^\omega$. Lemma 3.4 provides $\alpha_{0}\!\in\! G$ and 
$f\! :\! 2^\omega\!\rightarrow\! G$ continuous.\bigskip

\noindent $\bullet$ Let us show that $\rho^\xi_0(\alpha 
)\! =\!\rho^\xi_0\big( {\cal S}[\alpha_0\Delta f(\alpha )]\big)$, for each 
$1\!\leq\!\xi\! <\!\omega_1$ and $\alpha\!\in\! 2^\omega$. For $\xi\! =\! 1$ we 
apply Lemma 3.4.(b) to $t\!\in\!\omega$. Then we have, by induction:
$$\rho^{\eta +1}_0(\alpha )\! =\!\rho_0\big(\rho^{\eta }_0(\alpha )\big)\! =\! 
\rho_0\Big(\rho^{\eta }_0\big( {\cal S}[\alpha_0\Delta f(\alpha )]\big)\Big)\! =
\!\rho^{\eta +1}_0\big( {\cal S}[\alpha_0\Delta f(\alpha )]\big).$$
From this we deduce, by induction again, that
$$\rho^{(0,1)}_0(\alpha )\! =\!\rho^{\xi^\lambda_0}_0(\alpha )\! =\! 
\rho^{\xi^\lambda_0}_0\big( {\cal S}[\alpha_0\Delta f(\alpha )]\big)\! =\! 
\rho^{(0,1)}_0\big( {\cal S}[\alpha_0\Delta f(\alpha )]\big).$$
Thus $\rho^{(0,k+1)}_0(\alpha )\! =\!\rho^{(0,k+1)}_0\big( 
{\cal S}[\alpha_0\Delta f(\alpha )]\big)$, and
$$\rho^{\lambda}_0(\alpha )(k)\! =\! 
\rho^{(0,k+1)}_0(\alpha )(k)\! =\! 
\rho^{(0,k+1)}_0\big( {\cal S}[\alpha_0\Delta f(\alpha )]\big)(k)\! =\! 
\rho^{\lambda}_0\big( {\cal S}[\alpha_0\Delta f(\alpha )]\big)(k).$$ 
$\bullet$ This implies that $\alpha\!\in\! H_{1+\xi}$ is equivalent 
to ${\cal S}[\alpha_{0}\Delta f(\alpha )]\!\in\! H_{1+\xi}$ (for $\xi\!  
=\! 0$ we apply Lemma 3.4.(b) to $t\! :=\!\emptyset$).\bigskip 

\noindent $\bullet$ We argue by contradiction to show that 
$\lceil T\rceil\!\setminus\! S_{1+\xi}$ (resp., $S_{1+\xi}$) is not 
$\hbox{\rm pot}(\borapx )$ (resp., $\hbox{\rm pot}(\bormpx )$): there is a 
dense $G_{\delta}$ subset $G$ of $2^\omega$ such that 
$(\lceil T\rceil\!\setminus\! S_{1+\xi})\cap 
G^2$ (resp., $S_{1+\xi}\cap G^2$) is a $\borapx$ (resp., $\bormpx$) subset of 
$G^2$. But by the previous point we get $H_{1+\xi}\! =\! f^{-1}\big( [(\lceil 
T\rceil\!\setminus\! S_{1+\xi})\cap G^2]_{\alpha_{0}}\big)$ (resp., 
$\neg H_{1+\xi}\! =\! f^{-1}([S_{1+\xi}\cap G^2]_{\alpha_{0}})$), 
which is absurd.\hfill{$\square$}

\vfill\eject

\section{$\!\!\!\!\!\!$ Proof of Theorem 1.14.}\indent

 As announced in the introduction, we show more than Theorem 1.14.\bigskip

\noindent\bf Notation.\rm ~Let $X$, $Y$ be recursively presented 
Polish spaces. We set
$$B^{X\times Y}_0\! :=\!\{ p\!\in\! W^{X\times Y}\mid\exists (m,n)\!\in\! 
W^X\!\times\! W^Y~\ C^{X\times Y}_p\! =\! C^X_m\!\times\! C^Y_n\}.$$
Then we define an inductive operator $\Phi$ over $\omega$ (see [C]) 
as follows:\bigskip

\leftline{$\Phi (A)\! :=\! B^{X\times Y}_0\cup A\ \cup$}\bigskip

\rightline{$\{  p\!\in\! W^{X\times Y}\mid\exists\alpha\!\in\!\Borel\ 
\forall n\!\in\!\omega ~\ \ \alpha (n)\!\in\! W^{X\times Y}\cap A~\ 
\hbox{\rm and}~\ \neg C^{X\times Y}_p\! =\!\bigcup_n\ C^{X\times Y}_{\alpha (n)} 
\}.$}\smallskip\smallskip\noindent
Then $\Phi$ is clearly a $\Ca$ monotone inductive operator. We let, 
for any ordinal $\xi$, 
$${B^{X\times Y}_\xi\! =\!\Phi^\xi\! :=\!\Phi (\bigcup_{\eta <\xi}\ \Phi^\eta )}$$ 
(which is coherent with the definition of $B^{X\times Y}_0$).

\begin{thm} (Debs-Lecomte-Louveau) Let $T$ given by Theorem 1.11, 
$\xi\! <\!\omega^{\hbox{\it CK}}_{1}$, $S$ given by Theorem 1.11, and $X$, $Y$ be 
recursively presented Polish spaces.\medskip

\noindent $\bullet$ Let $A$, $B$ be disjoint $\Ana$ subsets of $X\!\times\! Y$. 
The following are equivalent:

\noindent (a) The set $A$ cannot be separated from $B$ by a 
$\hbox{\it pot}(\bormpx  )$ set.

\noindent (b) The set $A$ cannot be separated from $B$ by a $\Borel\cap 
\hbox{\it pot}(\bormpx )$ set.

\noindent (c) The set $A$ cannot be separated from $B$ by a $\bormpx (\tau_1)$ 
set.

\noindent (d) $\overline{A}^{\tau_{1+\xi}}\cap B\!\not=\!\emptyset$.

\noindent (e) There are $u\! :\! 2^\omega\!\rightarrow\! X$ and 
$v\! :\! 2^\omega\!\rightarrow\! Y$ continuous with 
$S\!\subseteq\! (u\!\times\! v)^{-1}(A)$ and 
$\lceil T\rceil\!\setminus\! S\!\subseteq\! (u\!\times\! v)^{-1}(B)$.\medskip

\noindent $\bullet$ The sets $W^{X\times Y}_{0}\! =\! B^{X\times Y}_{0}$, 
$W^{X\times Y}_{1+\xi}\! =\! B^{X\times Y}_{1+\xi}$ and $W^{X\times Y}_{<1+\xi}$ 
are $\Ca$.\end{thm}

\noindent\bf Proof.\rm ~The set $B^{X\times Y}_{0}$ is clearly $\Ca$ 
and a subset of $W^{X\times Y}_{0}$. Conversely, if $p$ is in 
$W^{X\times Y}_{0}$, then $C^{X\times Y}_{p}$ is a $\Ana$ rectangle, 
and a $\Borel$ rectangle by reflection. So $p\!\in\! B^{X\times 
Y}_{0}\! =\! W^{X\times Y}_{0}$.\bigskip

\noindent $\bullet$ We argue by induction on $\xi$. So assume that the result 
has been shown for $\eta\! <\!\xi$.\bigskip

\noindent $\bullet$  Let us show that $W^{X\times Y}_{<1+\xi}$ is $\Ca$. We may 
assume that $\xi\! =\! 1\! +\!\xi$ is an infinite limit ordinal since 
$W^{X\times Y}_{<\eta+1}\! =\! W^{X\times Y}_{\eta}$. By Lemma 4.8 in [C] 
the following relation is $\Ca$:
$$R(p,\delta )\ \Leftrightarrow\ \delta\!\in\!\hbox{\rm WO}\ \hbox{\rm and}\ 
p\!\in\!\Phi^{|\delta |}.$$
The following argument can be found in [Lo1], Proposition 1.4. Let 
$\delta_{\xi}\!\in\!\hbox{\rm WO}\cap\Borel$ with $|\delta_{\xi}|\! =\! 
\xi$, and $\delta_{\xi}^m$ be the restriction of the ordering 
$\delta_{\xi}$ to the $\delta_{\xi}$-predecessors of $m$. We get, by 
induction assumption,
$$\begin{array}{ll} 
p\!\in\! W^{X\times Y}_{<1+\xi}\  
& \!\!\!\!\Leftrightarrow\ \exists\eta\! <\!\xi\ \ ~p\!\in\! W^{X\times 
Y}_{\eta}\Leftrightarrow\ \exists\eta\! <\!\xi\ ~p\!\in\! B^{X\times Y}_{\eta}\cr &Ê\cr
& \!\!\!\!\Leftrightarrow\ \exists\eta\! <\!\xi\ ~\ p\!\in\! \Phi^{\eta}\ 
\Leftrightarrow\ \exists m\!\in\!\omega ~\ R(p,\delta_{\xi}^m).
\end{array}$$
This shows that $W^{X\times Y}_{<1+\xi}$ is $\Ca$.

\vfill\eject

\noindent (a) $\Rightarrow$ (b) and (a) $\Rightarrow$ (c) are clear since 
${\it\Delta}_X$ and ${\it\Delta}_Y$ are Polish.\bigskip

\noindent (c) $\Rightarrow$ (d) This comes from the fact that 
${\bormone(\tau_{1+\xi})\!\subseteq\!\bormpx (\tau_{1})}$.\bigskip

\noindent (d) $\Rightarrow$ (e) This comes from Theorems 2.4.1 and 2.4.4 (Lemma 
2.2.2 is at this moment true until the level $1\! +\!\xi$).\bigskip

\noindent (e) $\Rightarrow$ (a) If $D\!\in\!\hbox{\rm pot}(\bormpx )$ separates 
$A$ from $B$, then $S\! =\! (u\!\times\! v)^{-1}(D)\cap\lceil T\rceil$ is 
$\hbox{\rm pot}(\bormpx )$, which contradicts Theorem 1.11.\bigskip

\noindent (b) $\Rightarrow$ (d) We argue by contradiction, so that 
$\overline{A}^{\tau_{1+\xi}}$ separates $A$ from $B$. By induction 
assumption and the first reflection theorem there is $\alpha\!\in\!\Borel$ 
with $\alpha (n)\!\in\! W^{X\times Y}_{<1+\xi}$ and $C_{\alpha 
(n)}^{X\times Y}\!\subseteq\!\neg A$, for each integer $n$, and 
$B\!\subseteq\! E\! :=\!\bigcup_{n}\ C_{\alpha (n)}^{X\times Y}$. But $E$ is 
${\Borel\cap\hbox{\rm pot}(\borapx )}$ and separates $B$ from $A$, which is 
absurd.\bigskip

\noindent $\bullet$ The proof of the implication (b) $\Rightarrow$ (d) imply that 
$W^{X\times Y}_{1+\xi}$ is $\Ca$ since $W^{X\times Y}_{<1+\xi}$ is 
$\Ca$ and
$$W^{X\times Y}_{1+\xi}\! =\!\{ p\!\in\! W^{X\times Y}\mid 
\exists\alpha\!\in\!\Borel\ \forall n\!\in\!\omega\ ~\alpha 
(n)\!\in\! W^{X\times Y}_{<1+\xi}~\ \hbox{\rm and}~\ \neg C_p^{X\times 
Y}\! =\!\bigcup_{n}\ C_{\alpha (n)}^{X\times Y}\}.$$
$\bullet$ It remains to see that 
$W^{X\times Y}_{1+\xi}\! =\! B^{X\times Y}_{1+\xi}$. But by induction assumption 
we get
$$\begin{array}{ll}
& B^{X\times Y}_{1+\xi}\cr & \cr =\!\!\!\! 
& \Phi (\bigcup_{\eta <1+\xi}\ \Phi^\eta ) =\Phi (\bigcup_{\eta 
<1+\xi}\ B^{X\times Y}_{\eta})\cr & \cr =\!\!\!\! 
& \bigcup_{\eta <1+\xi}\ B^{X\times Y}_{\eta}\ \ \cup
\{ p\!\in\! W^{X\times Y}\mid\exists\alpha\!\in\!\Borel\ \forall 
n\!\in\!\omega\ ~\alpha (n)\!\in\!\bigcup_{\eta <1+\xi}\ B^{X\times 
Y}_{\eta}\ ~\hbox{\rm and}\cr 
& \ \ \ \ \ \ \ \ \ \ \ \ \ \ \ \ \ \ \ \ \ \ \ \ \ \ \ \ \ \ \ \ \ \ \ \ \ \ 
\ \ \ \ \ \ \ \ \ \ \ \ \ \ \ \ \ \ \ \ \ \ \ \ \ \ \ \ \ \ \ \ \ \ \ \ \ \ 
\ \ \ \ \ \ \ \ \ \ \ \ \ \ \ \ \ \ \ \ \ \ \ \ \ \ \ \ \ \ \ \ \ \ \ \ \ \ 
\neg C^{X\times Y}_p\! =\!\bigcup_n\ C^{X\times Y}_{\alpha (n)}\}\cr & \cr =\!\!\!\! 
& W^{X\times Y}_{<1+\xi}\ \ \cup
\{ p\!\in\! W^{X\times Y}\mid\exists\alpha\!\in\!\Borel\ \forall 
n\!\in\!\omega\ ~\alpha (n)\!\in\! W^{X\times Y}_{<1+\xi}\ ~\hbox{\rm and}~\ 
\neg C^{X\times Y}_p\! =\!\bigcup_n\ C^{X\times Y}_{\alpha (n)}\}\cr & \cr =\!\!\!\! 
& W^{X\times Y}_{1+\xi}.
\end{array}$$
This finishes the proof.\hfill{$\square$}\bigskip

\noindent\bf Remark.\rm\ As we saw with Theorem 2.2.1, the equivalence 
between (a), (b) and (c) is essentially shown in [Lo2]. It is also essentially 
shown in [Lo2] that (a), (b) and (c) are equivalent to (d) (see the 
proof of Theorem 2.8, (a) page 25, in [Lo2]). An immediate consequence of 
Theorem 4.1 is the following, shown in [Lo2]:

\begin{cor} (Louveau) Let $\xi\! <\!\omega^{\hbox{\it CK}}_{1}$, 
$X$, $Y$ be recursively presented Polish spaces, and $A$ a $\Borel$ subset of 
$X\!\times\! Y$. The following are equivalent:

\noindent (a) The set $A$ is $\hbox{\it pot}(\bormpx )$.

\noindent (b) The set $A$ is $\bormpx (\tau_1)$.\end{cor}

\vfill\eject

\section{$\!\!\!\!\!\!$ References.}

\noindent [B]\ \ B. Bollob\'as,~\it Modern graph theory,~\rm 
Springer-Verlag, New York, 1998

\noindent [C]\ \ D. Cenzer, Monotone inductive definitions over the continuum,~
\it J. Symbolic Logic\rm ~41 (1976), 188-198

\noindent [D-SR]\ \ G. Debs and J. Saint Raymond, Borel liftings of Borel sets: 
some decidable and undecidable statements,~\it to appear in Mem. Amer. Math. Soc.\rm

\noindent [H-K-Lo]\ \ L. A. Harrington, A. S. Kechris and A. Louveau, A Glimm-Effros 
dichotomy for Borel equivalence relations,~\it J. Amer. Math. Soc.\rm ~3 (1990), 903-928

\noindent [Hj-K-Lo]\ \ G. Hjorth, A. S. Kechris and A. Louveau, Borel equivalence 
relations induced by actions of the symmetric group,~\it Ann. Pure Appl. Logic\rm 
~92 (1998), 63-112

\noindent [K]\ \ A. S. Kechris,~\it Classical Descriptive Set Theory,~\rm 
Springer-Verlag, 1995

\noindent [L1]\ \ D. Lecomte, Classes de Wadge potentielles et 
th\'eor\`emes d'uniformisation partielle,\it ~Fund. Math.~\rm 143 (1993), 231-258

\noindent [L2]\ \ D. Lecomte, Uniformisations partielles et crit\`eres \`a la 
Hurewicz dans le plan,~\it Trans. Amer. Math. Soc.\rm ~347, 11 (1995), 4433-4460

\noindent [L3]\ \ D. Lecomte, Tests \`a la Hurewicz dans le plan,\it ~Fund. Math.~
\rm 156 (1998), 131-165

\noindent [L4]\ \ D. Lecomte, Complexit\'e des bor\'eliens~\`a coupes 
d\'enombrables,\ \it Fund. Math.~\rm 165 (2000), 139-174

\noindent [L5]\ \ D. Lecomte, On minimal non potentially closed subsets of the plane,
\ \it to appear in Topology Appl.~\rm 

\noindent [L6]\ \ D. Lecomte, Hurewicz-like tests for Borel subsets of the plane,
~\it Electron. Res. Announc. Amer. Math. Soc.\rm\ 11 (2005)

\noindent [Lo1]\ \ A. Louveau, A separation theorem for $\Ana$ sets,\ \it Trans. 
Amer. Math. Soc.\ \rm 260 (1980), 363-378

\noindent [Lo2]\ \ A. Louveau, Ensembles analytiques et bor\'eliens dans les 
espaces produit,~\it Ast\'erisque (S. M. F.)\ \rm 78 (1980)

\noindent [Lo-SR]\ \ A. Louveau and J. Saint Raymond, Borel classes and closed games: 
Wadge-type and Hurewicz-type results,\ \it Trans. Amer. Math. Soc.\ \rm 304 (1987), 431-467

\noindent [M]\ \ Y. N. Moschovakis,~\it Descriptive set theory,~\rm North-Holland, 1980

\noindent [SR]\ \ J. Saint Raymond, La structure bor\'elienne d'Effros 
est-elle standard ?,\ \it Fund. Math.\ \rm 100 (1978), 201-210

\end{document}